\documentclass[10pt]{amsart}

\usepackage{amsthm, amsmath}
\usepackage{amssymb}
\usepackage[all]{xy}
\numberwithin{equation}{subsection}
\DeclareMathOperator{\End}{End}
\DeclareMathOperator{\Lie}{Lie}
\DeclareMathOperator{\meas}{meas}
\DeclareMathOperator{\univ}{univ}

\DeclareMathOperator{\gen}{gen}

\theoremstyle{plain}

\newtheorem{theorem}{Theorem}[subsection]

\newtheorem{lemma}[theorem]{Lemma}

\newtheorem{cor}[theorem]{Corollary}

\newtheorem{rem}[theorem]{Remark}

\begin{document}

\title{Iwahori-Hecke algebras}

\author{Thomas J. Haines \and Robert E. Kottwitz \and Amritanshu Prasad}

\maketitle

Our aim here is to give a fairly self-contained
exposition of some basic facts about the Iwahori-Hecke
algebra $H$ of a split $p$-adic group $G$, including Bernstein's
presentation and description of the center, Macdonald's formula, the
Casselman-Shalika formula, and the Lusztig-Kato formula.

There are no new results here, and the same is
essentially true of the proofs. We have been strongly influenced
by the notes \cite{bernstein} of a course given by Bernstein.  In the spirit of Bernstein's work, we approach the material with an emphasis on the ``universal unramified principal series'' module $M = C_c(A_\mathcal O N\backslash G/I)$, which is a right module over the Iwahori-Hecke algebra $H = C_c(I\backslash G/I)$.  We use $M$ to develop the theory of intertwining operators in a purely algebraic framework.  Once this framework is established, we adapt it to produce rather efficient proofs of the above results, following closely at times earlier proofs.  In particular, in our treatment of Macdonald's formula and the Casselman-Shalika formula, we follow the method introduced by Casselman \cite{C'} and Casselman-Shalika \cite{CS}.  We follow Kato's strategy from \cite{kato} in proving a fundamental formula of Lusztig \cite{lusztig83}, which nowadays has come to be known as the Lusztig-Kato formula.  

The reader may find in \cite{Ram} another survey article which proves some of the results of the present paper by different methods.

\medskip

The following notation will be used throughout this paper.
We work over a $p$-adic field $F$ with valuation ring
$\mathcal O$ and prime ideal $P=(\pi)$. We denote by $k$ the residue field
$\mathcal O/P$ and by $q$ the cardinality of $k$.

Consider a split connected reductive group $G$ over $F$, with split
maximal torus $A$ and Borel subgroup $B=AN$ containing $A$. We
write $\bar B=A\bar N$ for the Borel subgroup containing $A$ that
is opposite to $B$. We assume that
$G,A,N$ are defined over~$\mathcal O$. We write
$K$ for $G(\mathcal O)$ and $I$ for the Iwahori subgroup of $K$
defined as the inverse image under
$G(\mathcal O) \to G(k)$ of $B(k)$. For $\mu \in X_*(A)$ we write
$\pi^\mu$ for the element
$\mu(\pi)
\in A(F)$. Note that $\mu \mapsto \pi^\mu$ gives an isomorphism
from $X_*(A)$ to $A/A_\mathcal O$. (We will often abbreviate $A(F)$
to $A$ and $A(\mathcal O)$ to $A_\mathcal O$, etc.)

\section{Bernstein's presentation \cite{lusztig89}}
\subsection{Extended affine Weyl group}
The extended affine Weyl group $\widetilde W$ is the quotient of
$N_{G(F)}(A)$ by $A_\mathcal O$. Thus $\widetilde W$ contains the
translation subgroup $A/A_\mathcal O=X_*(A)$, as well as the finite
Weyl group $W$, which we realize inside $\widetilde W$ as the quotient
of $N_K(A)$ by $A_\mathcal O$. Recall that $\widetilde W$ is the
semidirect product of $W$ and $X_*(A)$. In the early part of the paper, when
we are thinking about a cocharacter $\mu$ as an element of the translation
subgroup of~$\widetilde W$, we denote it by $\pi^\mu$; later on we often denote
it instead by $t_\mu$. 

\subsection{Iwahori-Hecke algebra H} We denote by $H$
 the Iwahori-Hecke algebra $C_c(I\backslash G/I)$. The convolution product
is defined using the Haar measure giving $I$ measure 1.  The elements
$T_x:=1_{IxI}$
($x \in \widetilde W$) form a
$\mathbf{C}$-basis for $H$. (Throughout the paper we write 
the characteristic
function of a subset $S$ as $1_S$.)

\subsection{Double cosets $A_\mathcal O N\backslash G/I=\widetilde W$}
The obvious map $\widetilde W
\to A_\mathcal O N\backslash G/I$ is a bijection.
How
does this decomposition
work? Write $g \in G$ as
$g=\pi^\mu n k$. Then write $k=n_\mathcal O wi$ with $n_\mathcal O
\in N_\mathcal O$, $i \in I$, $w
\in W$, with $w$ inside $K$. Thus $g=\pi^\mu nn_\mathcal O   w i$,
showing that
$g$ lies in the double coset of $\pi^\mu w \in \widetilde W$.  In short,
we read off
$\pi^\mu$ from the Iwasawa decomposition and then read off $w \in
W$ by applying the Bruhat decomposition over the residue field to
the element $k$ produced from the Iwasawa decomposition.

\subsection{Definition of the module $M$}
Put $M:=C_c(A_\mathcal O N\backslash
G/I)$. Note that $M$ is a right
$H$-module (since it arises as the $I$-fixed vectors in the smooth
$G$-module considered in \eqref{seconddescription} below).\footnote{
More explicitly, the action of  $h \in H$ on $m \in M$ is given by
the convolution $m \cdot h$ of the two functions, where convolution is
defined using the Haar measure on $G$  giving $I$ measure $1$. }
 For $x \in
\widetilde W$ we denote by $v_x$ the characteristic function
$1_{A_\mathcal O N x I}$. The elements $v_x$ ($x \in \widetilde W$)
form a $\mathbf C$-basis for $M$.  Of special importance (see Lemma
\ref{rank1} below) is the
basis element $v_1=1_{A_\mathcal O N  I}$.
Let $R=C_c(A/A_\mathcal O)=\mathbf C[X_*]$, the group algebra of
$X_*$, an abbreviation for
$X_*(A)$. Thus the elements $\pi^\mu$ ($\mu \in X_*$) form a basis
for the vector space $R$. We make $M$ into a left $R$-module as
follows. Let $\mu \in X_*$ and let $x \in \widetilde W$. Then
$\pi^\mu \cdot v_x:= q^{-\langle \rho,\mu \rangle} v_{\pi^\mu \cdot
x}$, where $\rho$ is half the sum of the roots of $A$ in
$\Lie(N)$.  (Note that the scalar
$q^{-\langle \rho,\mu \rangle}$ is equal to
$\delta_B(\pi^\mu)^{1/2}$, where for $a
\in A$, $\delta_B(a)$ denotes the absolute value of the determinant
of the adjoint action of $a$ on
$\Lie(N)$.) The actions of $R$ and $H$ commute, so that $M$ is an
$(R,H)$-bimodule.

\subsection{Second point of view on the module $M$} \label{sec.pt.view}
Consider the representation
(by right translations) of $G$ on $C^\infty_c(A_\mathcal O N\backslash G)$. 
It is compactly induced from the trivial 
representation of $A_\mathcal O N$; doing
the induction in stages, we see that
\begin{equation}\label{seconddescription}
C^\infty_c(A_\mathcal O N\backslash G)=i^G_B(R).
\end{equation}
Here we are using normalized induction and $R$ is viewed as $A$-module via
$\chi^{-1}_{\univ}$, where $\chi_{\univ}$ is the tautological character
$A/A_\mathcal O \to R^\times$ mapping $\pi^\mu$ to $\pi^\mu$.
 Thus an element
in this induced representation is a locally constant $R$-valued function $\phi$
on~$G$ satisfying
\begin{equation*}
\phi(ang)=\delta_B(a)^{1/2}\cdot a^{-1} \cdot \phi(g)
\end{equation*}
for all $a \in A$, $n \in N$, $g \in G$, and the group $G$ acts by right
translations.
The isomorphism \eqref{seconddescription} has the following explicit
description. Let $\varphi \in C^\infty_c(A_\mathcal O N\backslash G)$. Then the
corresponding element $\phi \in i^G_B(R)$ is defined by
\begin{equation}
\phi(g)=\sum_{a \in A/A_\mathcal O} \delta_B(a)^{-1/2} \varphi(ag) \cdot a
\end{equation}
for $g \in G$. There is an obvious $R$-module structure on $i^G_B(R)$, with
$r\phi$ given by $(r\phi)(g)=r( \phi(g))$. The isomorphism
\eqref{seconddescription} induces an $(R,H)$-bimodule isomorphism from $M$ to
the Iwahori fixed vectors in $i^G_B(R)$.

Let $\chi$ denote a quasicharacter $\chi:A/A_\mathcal O \to \mathbf
C^\times$. Then $\chi$ determines a $\mathbf C$-algebra
homomorphism $R \to \mathbf C$. Using this homomorphism to extend
scalars, we obtain the $H$-module
\begin{equation*}
\mathbf C \otimes_R M = \mathbf C \otimes_R i^G_B(\chi^{-1}_{\univ})^I
=i^G_B(\chi^{-1})^I.
\end{equation*}

\subsection{Structure of the module $M$}

The next result  is due to Chriss and
Khuri-Makdisi \cite{chriss}, who derived it from 
Bernstein's presentation of the
Iwahori-Hecke algebra. Here we turn the logic around, first studying $M$
directly, then using it to produce Bernstein's presentation.

\begin{lemma}\label{rank1}
The map $h \mapsto v_1 h$ is an
isomorphism of right $H$-modules from $H$ to
$M$. In other words $M$ is free of rank $1$ as $H$-module with
canonical generator $v_1$. In particular we have a canonical
isomorphism $H \simeq \End_H(M)$, which identifies $h' \in H$ with
the endomorphism $v_1h \mapsto v_1h'h$  of $M$.
\end{lemma}

\begin{proof}  It suffices to show that the map $h \mapsto v_1h$,
written in terms of the bases $\{ T_w \}_w$ and $\{ v_w \}_w$, is
``a triangular matrix with non-zero diagonal''.  This follows from the
following claim.

\noindent{\bf Claim:} $NxI \cap IyI \neq \emptyset  \Rightarrow x
\leq y \,\,\, \mbox{in the Bruhat ordering}$.\footnote{The Iwahori
subgroup $I$ determines in a canonical way the Bruhat order on $\tilde
W$, but only when $\widetilde W$ is viewed as $N_G(A)/A_\mathcal O$. When  
$\widetilde W$ is viewed as the semidirect product of $W$ and $X_*(A)$,  the
Bruhat order depends on the normalization of the isomorphism between
$X_*(A)$ and
$A/A_\mathcal O$. Our normalization is 
$\mu \mapsto \pi^\mu$, and therefore our Bruhat order on
$X_*(A) \rtimes W$  is the one determined by the simple affine
reflections about the walls of the unique alcove in $X_*(A)$ whose
closure contains the origin and lies in the \emph{negative} Weyl 
chamber.  See section \ref{Kato}.}

{\em Proof of
Claim:}  Suppose $nx \in IyI$, for $n \in N$.  Choose $\mu$ so
dominant that $\pi^\mu n \pi^{-\mu} \in I$.  Then $(\pi^\mu n
\pi^{-\mu}) \pi^\mu x \in \pi^\mu I yI$, hence
$$ I\pi^\mu x I \subset I\pi^\mu I y I \subset \coprod_{y' \leq y}
I \pi^\mu y' I,
$$ from which the claim follows.

\end{proof}

The following three equalities (see \cite{chriss}, \cite{reeder}) are also
useful.
\begin{align}
v_1 T_w &= v_w, \text{ for every $w \in W$,} \label{use1} \\
v_{\pi^\mu} T_w &= v_{\pi^\mu w}, \text{ for every $w \in W$ and
$\mu \in X_*(A)$,} \label{use45} \\
v_1 T_{\pi^\mu} &= v_{\pi^\mu}, \text{ for $\mu \in X_*(A)$
dominant.} \label{use2}
\end{align}
Recall  the
Iwahori factorization $I=(I \cap N)A_\mathcal O (I \cap \bar N)$.
The first equality uses $A_\mathcal O NI \cdot IwI=A_\mathcal O NwI$ (a
consequence of the Iwahori factorization) as well as $A_\mathcal O NI \cap
wIw^{-1}I =I $ (a consequence of $A_\mathcal O NI \cap
K =I $),
and  the second equality follows from the first (using the left $R$-module
structure on $M$).  The  third equality uses $A_\mathcal O NI
\cdot I\pi^\mu I=A_\mathcal O N\pi^\mu I$, a
consequence of the Iwahori factorization and the
dominance of $\mu$,  which implies that
\[
\pi^\mu(I \cap N)\pi^{-\mu} \subset I\cap N
\,\text{ and }\,
\pi^{-\mu}(I \cap \bar N)\pi^{\mu} \subset I\cap \bar N,
\]
and also uses $A_\mathcal O NI \cap
\pi^\mu I\pi^ {-\mu}I =I $, which we leave as an exercise for the reader.

\subsection{Rough structure of the algebra $H$}
The finite dimensional Hecke algebra $H_0=C(I\backslash K/I)$ is a
subalgebra of $H$. Moreover, elements in
$R$ can be viewed as endomorphisms of $M$, and hence by the
previous lemma can be considered as elements in $H$. In this way we
embed $R$ as a subalgebra of $H$. We will denote by $\Theta_\lambda \in H$ the
image of the basis element $\pi^\lambda$ of~$R$ under the embedding $R
\hookrightarrow H$. Unwinding the definitions, one finds the basic identity
\begin{equation}
v_1\Theta_\lambda=\pi^\lambda  v_1,
\end{equation}
which says that $v_1$ is an eigenvector for the right action of the subalgebra
$R$ of~$H$.

\begin{lemma}\label{bases}
 Multiplication in $H$ induces a vector space isomorphism
\[
R \otimes_\mathbf{C} H_0\overset{\approx} \longrightarrow H,
\]
 sending $\pi^\mu \otimes h$
to
$\Theta_\mu h$.  Composing this isomorphism with the isomorphism $h \mapsto
v_1 h$ considered above, we get a vector space isomorphism from $R
\otimes_\mathbf{C} H_0$ to $M$, sending $\pi^\mu \otimes T_w$ to $q^{-\langle
\rho, \mu \rangle} v_{\pi^\mu w}$.
\end{lemma}

\begin{proof}  Using \eqref{use1} and the definitions, one checks that the
composition
\[
R \otimes_\mathbf C H_0 \to H \to M
\]
sends
$\pi^\mu \otimes T_w$ to $q^{-\langle \rho, \mu \rangle}
v_{\pi^\mu w}$ and is hence an isomorphism. Since $H \to M$ is an
isomorphism by Lemma \ref{rank1}, the map $R \otimes_\mathbf C H_0 \to H$ is
also an isomorphism.

\end{proof}

\begin{rem}\label{rem172}
   It follows from \eqref{use2} that
$\Theta_\lambda$ agrees with the element denoted by this symbol in
Lusztig's work: namely,
$\Theta_\lambda = q^{\langle \rho, -\lambda_1 + \lambda_2 \rangle}
T_{\pi^{\lambda_1}} T^{-1}_{\pi^{\lambda_2}}$, where $\lambda =
\lambda_1 - \lambda_2$, and $\lambda_1$, $\lambda_2$ are dominant
cocharacters.
\end{rem}

\subsection{Involutions on $R$ and $H$}\label{invdefs}
Recall that in order to pass back and forth between left and right $G$-modules
one uses the anti-isomorphism $g \mapsto g^{-1}$ from $G$ to itself. The
corresponding way of passing from left to right $H$-modules uses the standard
anti-involution $\iota$ on~$H$ given by
$\iota (h)(x) = h(x^{-1})$.

Moreover there is also an involution
$\iota_A$ on $R$ (which is the Iwahori-Hecke algebra for $A$); thus
$\iota_A$ sends $\pi^\mu$ to $\pi^{-\mu}$.

\subsection{A sesquilinear form on $M$} There is an $R$-valued
pairing on $i^G_B(\chi^{-1}_{\univ})$, defined by
\begin{equation}
(\phi_1,\phi_2):=\oint_{B\backslash G} \iota_A(\phi_1(g))\phi_2(g).
\end{equation}
 What is the meaning of $\oint_{B\backslash G}$? Consider the induced
representation $i^G_B(\delta_B^{1/2})$, which consists of locally constant
functions $F$ on $G$ satisfying
\begin{equation*}
F(ang)=\delta_B(a)F(g).
\end{equation*}
The space of $G$-invariant linear functionals on $i^G_B(\delta_B^{1/2})$ is
$1$-dimensional; we denote by $\oint_{B\backslash G}$ the unique such functional
that takes the value $1$ on the function $F_0 \in i^G_B(\delta_B^{1/2})$ defined
by
$F_0(ank)=\delta_B(a)$.

This pairing is sesquilinear, in the sense that
\begin{equation}\label{sesq}
(r_1\phi_1,r_2\phi_2)=\iota_A(r_1)r_2\cdot (\phi_1,\phi_2).
\end{equation}
Moreover it satisfies
\begin{equation}\label{herm}
(\phi_2,\phi_1)=\iota_A(\phi_1,\phi_2)
\end{equation}
and is $G$-invariant.

Note that $\phi \mapsto \iota_A \circ \phi$ is an $\iota_A$-linear isomorphism
from $i^G_B(\chi^{-1}_{\univ})$ to $i^G_B(\chi_{\univ})$. Therefore our
sesquilinear form can also be thought of as an $R$-bilinear pairing
\begin{equation}\label{rbil}
i^G_B(\chi_{\univ}) \otimes_R i^G_B(\chi^{-1}_{\univ}) \to R.
\end{equation}

After extending scalars
$R \to \mathbf C$ using a quasicharacter
$\chi: A/A_\mathcal O \to
\mathbf C^\times$, the pairing \eqref{rbil} becomes the standard pairing
\begin{equation}
i^G_B(\chi) \otimes_\mathbf C i^G_B(\chi^{-1}) \to \mathbf C.
\end{equation}

Recall that we have identified $M$ with the Iwahori fixed vectors in
$i^G_B(\chi^{-1}_{\univ})$. Thus, by restriction, we get a perfect sesquilinear
form on $M$, which we denote by $(m_1,m_2)$. It satisfies the Hecke algebra
analog of $G$-invariance, namely
\begin{equation}\label{h-inv}
(m_1h,m_2)=(m_1,m_2\iota(h))
\end{equation}
for all $h \in H$.

\subsection{Generalities on intertwiners}\label{gener}
For each $w \in W$ we would
like to define an intertwiner $I_w$  from one suitable completion
of $M$ to another.  To this end it is best to let the Borel
subgroup vary (and then recover $I_w$ by bringing the second Borel
back to the first by an element of the Weyl group).  In this
discussion the maximal torus $A$, the Iwahori subgroup $I$, and the
maximal compact $K$ will remain fixed.  We let $\mathcal B(A)$
denote the set of Borel subgroups containing $A$.  For $B = AN \in
\mathcal B(A)$, put $M_B = C_c(A_{\mathcal O} N \backslash G/I)$.

 First let us discuss the
completions that will come up. Let $J$ be a set of coroots which is
a  subset of some system of positive coroots. As usual $R$ denotes
the group algebra of $X_*(A)$.  We denote by $\mathbf C[J]$ the
$\mathbf C$-subalgebra of $R$ generated by~$J$, and by $\mathbf
C[J]{\hat{\,}}$  the completion of $\mathbf C[J]$ with respect to
the (maximal) ideal generated by $J$. Finally, we denote by $R_J$
the $R$-algebra $\mathbf C[J]{\hat{\,}} \otimes_{\mathbf C[J]} R$,
a completion of~$R$ that can be viewed as the convolution algebra
of complex valued functions on~$X_*(A)$ supported on a finite union
of sets of the form $x+C_J$, with  $x\in X_*(A)$ and where
$C_J$ is the submonoid of $X_*(A)$ consisting of all non-negative
integral linear combinations of elements in~$J$.

  Given $B=AN \in
\mathcal B(A)$ and given $J$ as above, we then denote by $M_{B,J}$
the module $R_J \otimes_R M_B$, which  can be thought of as
consisting of functions $f$ on $A_\mathcal O N\backslash G/I$
satisfying the following support condition: there exists a finite
union $S$ of sets of the form $x + C_J$ such that the support of $f$ 
is contained in the union of the sets
$A_{\mathcal O}N \pi^\nu K$ for $\nu \in S$.  It is clear that
$M_{B,J}$ is a left $R_J$-module and a right $H$-module.

 Now let
$B=AN$,
$B'=AN'$ be two Borel subgroups in $\mathcal B(A)$. As usual we
write $\bar B=A\bar N$ for the Borel subgroup in $\mathcal B(A)$
opposite to~$B$. Let $J$ be the set of coroots that are positive
for $B'$ and negative for $B$. We are going to define an
intertwiner $I_{B',B}:M_{B,J}\to M_{B',J}$. This intertwiner is an
$(R_J,H)$-bimodule map, and is defined as follows (viewing elements
in completions as functions, as above). Let $\varphi \in M_{B,J}$.
Then the intertwiner $I_{B',B}$ takes $\varphi$ to the function
$\varphi'$ on $A_\mathcal O N'\backslash G/I$ whose value at $g \in
G$ is defined by the integral
\begin{equation*}
\varphi'(g)=\int_{N'\cap \bar N} \varphi (n'g) dn'.
\end{equation*}
The Haar measure $dn'$ is normalized to give $N' \cap \bar N \cap K$
measure~$1$.  Note that the integral makes sense since the
integrand is a smooth and compactly supported function on the group
$N'\cap \bar N$ (smoothness being trivial, compact support
requiring justification, to be done in the lemma below). In fact
things still work fine if we enlarge $J$ in any way (but so that
the enlarged set is still contained in some positive system, for
instance, the positive system defined by $B'$).

 Now suppose that
we have three Borel subgroups $B_1=AN_1$, $B_2=AN_2$, $B_3=AN_3$ in
$\mathcal B(A)$. Let $J_{ij}$ be the set of coroots that are
positive for~$B_i$ and negative for~$B_j$,  and assume that
$J_{31}$ is the disjoint union of $J_{21}$ and $J_{32}$. Write
$I_{ij}$ as an  abbreviation for the intertwiner $I_{B_i,B_j}$.
Then $I_{21}$, $I_{32}$ and $I_{31}$ can all  be defined using the
biggest of the three sets $J_{ij}$, namely $J_{31}$, and when this
is done we have the equality
\begin{equation*} I_{31}=I_{32}I_{21}.
\end{equation*}  In this formula we could also have taken $J$ to be
the set of all coroots that are positive for~$B_3$.

Why do the integrals make sense? For this we need the following
lemma, in which we return to $B$, $B'$ as above.

\begin{lemma}\label{compact}
For $\nu \in X_*(A)$ define a subset $C_\nu$ of the
group $N' \cap \bar N$ by $C_\nu := N' \cap \bar N
\cap\pi^{\nu}NK$.  Then:
\begin{enumerate}
\item If $C_\nu$ is non-empty, then $\nu$ is a non-negative
integral linear combination of coroots that are positive for~$B$
and negative for~$B'$.
\item The subset $C_\nu$ is  compact.
\end{enumerate}
\end{lemma}

\begin{proof}  We begin by recalling the definition of the
retraction $r_B:G\to X_*(A)$. Let $g \in G$ and use the Iwasawa
decomposition to write $g=\pi^\mu nk$ for some $\mu \in X_*(A)$,
$n\in N$, $k \in K$; then put $r_B(g):=\mu$. It is well-known that
$r_{B'}(g)-r_B(g)$ is a non-negative integral linear combination of
coroots that are positive for $B'$ and negative for~$B$.  (It is
enough to prove this for adjacent $B$, $B'$, for which a simple
computation in $SL(2)$ does the job.)

To prove the first statement
we consider an element $g \in C_\nu$. It is clear from the
definition of $C_\nu$ that $r_{B'}(g)=0$ and $r_B(g)=\nu$.
Therefore $\nu=r_B(g)-r_{B'}(g)$ is a non-negative integral linear
combination of coroots that are positive for $B$ and negative
for~$B'$.

 Now we turn to the second statement. It is enough to
prove that $\bar N \cap NC$ is compact for any compact subset $C$
of~$G$, which is equivalent to proving that the map $\bar N \to
N\backslash G$ is proper (in the topological sense). But in fact
$\bar N \to N\backslash G$ is a closed immersion (in the algebraic
sense), as follows from the fact that $N\bar N$ is closed in~$G$.
Recall the proof of this: For any dominant weight $\lambda$ there
exists a unique regular function $f_\lambda$ on the algebraic
variety $G$ such that $f_\lambda(na\bar n)=\lambda(a)^{-1}$ for all
$n \in N$, $a \in A$, $\bar n \in \bar N$. Then $N\bar N$ is the
closed subvariety defined by the equations $f_\lambda=1$ (one for
every dominant $\lambda$).
\end{proof}

Next we need to understand how the intertwiners behave with respect to the
sesquilinear form on $M_B$. Denote by $-J$ the set of negatives of the coroots
in~$J$. The involution $\iota_A$ on~$R$ extends to an isomorphism, still denoted
$\iota_A$, between $R_J$ and $R_{-J}$, and the sesquilinear form $(\cdot,\cdot)$
on~$M_B$ extends to our completions in the following sense: given $m_1 \in
M_{B,-J}$ and $m_2 \in M_{B,J}$ our old definition of $(m_1,m_2)$  still makes
sense and yields an element of $R_{J}$. The extended form $(\cdot,\cdot)$ still
satisfies
\eqref{sesq}.

Consider the
intertwiner $I_{B',B}:M_{B,J}\to M_{B',J}$, where $J$ denotes (as before)
 the set of coroots that are positive for $B'$ and negative for $B$.
 We also have the
intertwiner $I_{B,B'}:M_{B',-J}\to M_{B,-J}$.
 Let $m \in M_{B,J}$ and $m' \in M_{B',-J}$. Then we claim that
\begin{equation}\label{intadj}
(m',I_{B',B}m)=(I_{B,B'}m',m).
\end{equation}
Indeed, let $\phi$, $\phi'$ be the elements of $i^G_B(\chi^{-1}_{\univ})
\otimes_R
R_{J}$,
$i^G_{B'}(\chi_{\univ})\otimes_R R_{J}$ corresponding to $m$, $m'$
respectively.  Put $H:=A(N \cap N')$.  Then one sees
easily that both sides of the last equality are equal to
\begin{equation}
\oint_{H\backslash G} \phi'(g)\phi(g),
\end{equation}
where $\oint _{H\backslash G}$ is the unique $G$-invariant linear
functional on
\[
\{f \in C^\infty(G) :  f(hg)=\delta_{H}(h)f(g) \quad(\forall h \in H), \text{
compactly supported mod $H$} \}
\]
that takes the value $1$ on the function $f_0$ supported on $HK$ whose values
on~$HK$ are given  by
$f_0(hk)=\delta_{H}(h)$.

\subsection{Intertwiners $I_w$} We return now to the earlier
notation, where $B = AN$ is a fixed Borel subgroup. For each $w \in
W$, we define an intertwiner
\[
I_w: M_{B,w^{-1}J} \rightarrow M_{B,J}
\]
as the composition $I_{B,wB} L(w)$.  Here $L(w)$ is the
isomorphism
$M_{B,w^{-1}J} ~\tilde{\to}~ M_{wB,J}$ given by
$(L(w)\phi)(g) = \phi(\dot{w}^{-1}g)$,
where $\dot{w}$ is a
representative for $w$ taken in $K$.
Thus $I_w$ is defined by
the integral
\begin{equation}
I_w(\varphi)(g)=\int_{N_w} \varphi(\dot{w}^{-1}ng) \, dn,
\end{equation}
where $N_w$ denotes $N\cap w\bar N w^{-1}$.

{}From the discussion above,
the following properties are immediate.
\begin{lemma}\label{intprops}
We have
\begin{enumerate}
\item [(i)] $I_w \circ \pi^\mu = \pi^{w\mu} \circ I_w$, ~$\forall
\mu \in X_*(A)$,
\item [(ii)] $I_{w_1w_2} = I_{w_1} \circ I_{w_2}$,  if $l(w_1w_2) =
l(w_1) + l(w_2)$,
\item [(iii)] $I_w$ is a right $H$-module homomorphism.
\end{enumerate}
\end{lemma}

\subsection{Intertwiners in the rank $1$ case}

We suppose for the moment that $G$ has semisimple rank $1$. We
write $\alpha$ for the unique positive root of $A$, and $s_\alpha$
for the corresponding simple reflection, in this case the unique
non-trivial element in $W$.

Now we compute $\varphi'=I_{s_\alpha}(\varphi)$ for
$\varphi=v_1=1_{A_\mathcal O NI}$.  We write
$J(j,w)$ ($j
\in \mathbf Z$, $w \in W$) for the value of $\varphi'$ at the
element $\pi^{j\alpha^\vee}w$. Note that other values of $\varphi$
are $0$ and also that
$J(j,w)=0$ unless $j\ge 0$, which we now assume. At this point we
may as well take $G=SL(2)$. To simplify notation we temporarily
write $\mu$ for $j\alpha^\vee$.

  First suppose that $j=0$.  Note
that $s_\alpha nw \in A_\mathcal O NK$ iff
$n \in N_\mathcal O$. For $n \in N_\mathcal O$ the element
$s_\alpha nw$ belongs to $K$ and  hence belongs to
$A_\mathcal O NI$ iff its lower left entry is in the prime ideal in
$\mathcal O$. We conclude that
$J(0,1)=0$ and that $J(0,s_\alpha)=q^{-1}$.

Suppose $j > 0$. We have $s_\alpha=\begin{bmatrix} 0 & -1\\1 & 0
\end{bmatrix}$,
$n=\begin{bmatrix} 1 & x\\0 & 1 \end{bmatrix}$,
$\pi^\mu=\begin{bmatrix} \pi^j & 0\\0 &
\pi^{-j}\end{bmatrix}$,  so that
\begin{equation*} s_\alpha n \pi^{\mu} =
\begin{bmatrix} 0 & -\pi^{-j}\\\pi^j & x\pi^{-j} \end{bmatrix}.
\end{equation*} For $s_\alpha n \pi^{\mu}w$ to lie in $A_\mathcal O
NK$, we must have $x \in \pi^j\mathcal O^\times$. We now assume
this and write $x=\pi^ju$ for some unit $u$. Then
$s_\alpha n \pi^{\mu}=\begin{bmatrix} u^{-1} & -\pi^{-j}\\0 & u
\end{bmatrix}\begin{bmatrix}  1 & 0\\u^{-1} \pi^j & 1
\end{bmatrix}$, the first factor lying in $A_\mathcal O N$, the
second factor lying in $K$. Therefore $s_\alpha n \pi^{\mu} \in A_O
NI $ iff the second factor lies in $I$, which is always the case.
Therefore  $J(j,1)$ is the measure of
$\pi^j\mathcal O^\times$, namely $q^{-j}(1-q^{-1})$.

Moreover $s_\alpha n \pi^{\mu}s_\alpha \in A_\mathcal O NI$ iff
the  product of the second factor and $s_\alpha$, namely
$\begin{bmatrix}  0 & -1\\ 1 & -u^{-1} \pi^j  \end{bmatrix}$, lies
in $I$, which never happens. Therefore
$J(j,s_\alpha)=0$.

We have proved:
\begin{lemma}
 $\varphi'=q^{-1}v_{s_\alpha}+(1-q^{-1})\sum_{j=1}^\infty
q^{-j}v_{\pi^{j\alpha^\vee}}$.
\end{lemma}

Even easier:
\begin{lemma}
The intertwiner sends $1_{A_\mathcal O NK}$ to
 $$q^{-1}1_{A_\mathcal O NK}+\sum_{j=0}^\infty
q^{-j}(1-q^{-1})1_{A_\mathcal O
N\pi^{j\alpha^\vee}K}=\frac{1-q^{-1}\pi^{\alpha^\vee}}
{1-\pi^{\alpha^\vee}}1_{A_\mathcal O NK}.$$
\end{lemma}

\subsection{Consequences of the calculations above} Now we return to the general
case. In the next  lemma the calculations reduce easily to the rank $1$
case treated above, so we just record the results.

\begin{lemma}\label{lemma4.1} Let $\alpha$ be a simple root and $s_\alpha$ the
corresponding simple reflection. Then
\begin{gather}
I_{s_\alpha}(v_1)=q^{-1}v_{s_\alpha}+(1-q^{-1})\sum_{j=1}^\infty
{\pi^{j\alpha^\vee}}v_1. \tag{i} \\
I_{s_\alpha}(v_1+v_{s_\alpha})=\Bigl(\frac{1-q^{-1}\pi^{\alpha^\vee}}
{1-\pi^{\alpha^\vee}} \Bigr) (v_1+v_{s_\alpha}). \tag{ii} \\
I_{s_\alpha}(1_{A_\mathcal O NK})=\Bigl(\frac{1-q^{-1}\pi^{\alpha^\vee}}
{1-\pi^{\alpha^\vee}} \Bigr) 1_{A_\mathcal O NK}. \tag{iii}
\end{gather}

\end{lemma}

We now introduce the following notation. For $w \in W$ we denote by $R_w$ the
 set of positive roots $\alpha$ such that $w^{-1}\alpha$ is negative.

\begin{cor}[Gindikin-Karpelevich formula]\label{lemma4.2}
For $w \in W$ we have
\begin{equation*}
I_w(1_{A_\mathcal O NK})=\Bigl(\prod_{\alpha \in R_w} \frac{1-q^{-1}\pi^{\alpha^\vee}}
{1-\pi^{\alpha^\vee}} \Bigr) 1_{A_\mathcal O NK}.
\end{equation*}
\end{cor}

\subsection{Intertwiners $J_w$ without denominators} To
eliminate denominators we define a new intertwiner $J_w$ ($w \in
W$) by
$J_w:=\bigl( \prod_{\alpha \in
R_w}(1-\pi^{\alpha^{\vee}})\bigr) \cdot I_w$. Note
that $J_w$ preserves the subspace $M$ of $M_{B,w^{-1}J}$ and
$M_{B,J}$ and hence can be regarded as an element of $H$, via our
identification of $H$ with $\End_H(M)$. For a simple root $\alpha$,
the element of $H$ corresponding to $J_{s_\alpha}$ is (by Lemma
\ref{lemma4.1}(i))
equal to
\begin{equation} \label{589}
(1-q^{-1})\pi^{\alpha^\vee}+q^{-1}(1-\pi^{\alpha^{\vee}})T_{s_\alpha}
\footnote{Here we abuse notation and write $\pi^\lambda$ in place of its image $\Theta_\lambda$ under our embedding $R \hookrightarrow H$.  We will often do this (e.g. in (1.15.2) and again in section 2.1), leaving context to dictate what is really meant by $\pi^\lambda$.}.
\end{equation}

\subsection{Bernstein's relation}
Equation \eqref{589}, together with the equality
\begin{equation}\label{com-rel} J_w \circ \pi^{\mu} =\pi^{w(\mu)}
\circ J_w
\end{equation}  (for
$w = s_\alpha$), yields Bernstein's relation:
 \begin{equation}\label{bernrel}
 T_{s_\alpha}\pi^\mu=\pi^{s_\alpha(\mu)}
T_{s_\alpha}+(1-q)\frac {\pi^{s_\alpha(\mu)}-\pi^{\mu}}
{1-\pi^{-\alpha^\vee}}.
\end{equation}
 Using Bernstein's relation one can calculate the
square of $J_{s_\alpha}$, viewed as element in
$H$; it turns out to be the element
$(1-q^{-1}\pi^{\alpha^\vee})(1-q^{-1}\pi^{-\alpha^\vee})$ in the subalgebra $R$
 of~$H$.

Lemma \ref{bases} together with Bernstein's relation \eqref{bernrel} gives
Bernstein's presentation of~$H$.

\section{The center of $H$}
\subsection{A preliminary result}
 We are
going to prove that the subalgebra $R^W$ is the center of $H$, but
we start by proving something weaker.
\begin{lemma} The subalgebra $R^W$ is contained in the center of
$H$.
\end{lemma}
\begin{proof}
  Let $r \in R^W$. Then $r$ commutes with all elements in $R$, so
by Lemma \ref{bases} it suffices to show it commutes with
$T_{s_\alpha}$ for all simple
$\alpha$. By the intertwining property \eqref{com-rel} of
$J_{s_\alpha}$, it does commute with
$(1-q^{-1})\pi^{\alpha^\vee}+q^{-1}(1-\pi^{\alpha^{\vee}})T_{s_\alpha}$.
 So $r$
commutes with
$(1-\pi^{\alpha^{\vee}})T_{s_\alpha}$ and hence the bracket of $r$
and $T_{s_\alpha}$ is annihilated by $1-\pi^{\alpha^{\vee}}$.
Since $H$ is a free $R$-module,  the bracket
vanishes.
\end{proof}

\subsection{The normalized intertwiners $K_w$}
Let $L$ denote the field of fractions of the integral domain $R$.
Then $L^W$ is the field of fractions of $R^W$. We now consider the algebra
$H_{\gen}:=L^W\otimes _{R^W} H$ and the module $M_{\gen}:=L \otimes_R M=L^W
\otimes_{R^W} M$, which is an $(L,H_{\gen})$-bimodule.

 We define  the normalized
intertwiners by
\begin{equation}
K_w:=\Bigl(\prod_{\alpha \in R_w}\frac{1}{1-q^{-1}\pi^{\alpha^\vee}}\Bigr) \cdot
J_w=\Bigl(\prod_{\alpha \in
R_w}\frac{1-\pi^{\alpha^\vee}}{1-q^{-1}\pi^{\alpha^\vee}}\Bigr)
\cdot I_w.
\end{equation}
Each $K_w$ is an endomorphism of the $H_{\gen}$-module $M_{\gen}$ and fixes the
spherical vector $1_{A_\mathcal O NK}$, as one sees from Corollary
\ref{lemma4.2}. For simple
$\alpha$ we have $K_{s_\alpha}^2=1$. It follows from this and Lemma
\ref{intprops} that
\begin{equation}
K_{w_1w_2}=K_{w_1}K_{w_2}
\end{equation}
for all $w_1,w_2 \in W$.

The involution $\iota_A$ extends to~$L$, and
our sesquilinear pairing form $(\cdot,\cdot)$ on~$M$ extends to a sesquilinear
$L$-valued form, still denoted $(\cdot,\cdot)$, on ~$M_{\gen}$. It follows from
\eqref{intadj} that
\begin{equation}\label{223abc}
w\bigl(K_{w^{-1}}(m),m'\bigr)=\bigl(m,K_w(m')\bigr)
\end{equation}
for all $m,m' \in M_{\gen}$.

For later use we remark that
it follows from \eqref{589} that for any $w \in W$ one has
\begin{equation} \label{triang}
K_w(v_1)=\sum_{w' \le w} a_{ww'} \cdot  v_{w'}
\end{equation}
for certain elements  $a_{ww'} \in L$, with the diagonal elements
given by the simple formula
\begin{equation}\label{diag.elts}
a_{ww}=
\prod_{\alpha \in
R_w}\frac{1-\pi^{-\alpha^\vee}}{1-q\pi^{-\alpha^\vee}}.
\end{equation}

\subsection{Calculation of the center of $H$}
Since the endomorphism ring of the $H_{\gen}$-module $M_{\gen}$ is $H_{\gen}$,
we can view the endomorphisms $K_w$ as elements of $H_{\gen}$.
  The map $w \mapsto K_w$ is a group
homomorphism from $W$ to $H_{\gen}^\times$ and therefore induces an
algebra homomorphism from the twisted group algebra $L[W]$ to $H_{\gen}$.
\begin{lemma} \label{center}
The homomorphism $L[W] \to H_{\gen}$ is an isomorphism. The
center of $H_{\gen}$ is $L^W$. The center of $H$ is $R^W$.
\end{lemma} 
\begin{proof} The twisted group algebra is a matrix algebra over
$L^W$, and is therefore simple, which implies our map is
injective.   Comparing dimensions, we see that the map is an
isomorphism. Therefore $H_{\gen}$ is a matrix algebra over $L^W$, and
its center is
$L^W$. It follows easily that the center of $H$ is $R^W$. (Use
along the way the obvious fact that
$H$ is torsion-free as
$R^W$-module.)
\end{proof}

\section{Application: Restriction of two involutions to the center}
\subsection{Restriction of $\iota$ to the center}
Recall from before the anti-involution $\iota:H \to H$ given by
$\iota (h)(x) = h(x^{-1})$.  We are going to see that the
restriction of $\iota$ to the center of $H$ is very simple.

\begin{lemma}\label{iota}
 There are two involutions on $R^W$, one obtained by
restricting $\iota_A$ to $R^W$, the other obtained by restricting
$\iota$ to the center of $H$, which we have identified with $R^W$.
The  two involutions on $R^W$ coincide.
\end{lemma}

\begin{proof}  This follows from \eqref{sesq}, \eqref{h-inv}, the
non-degeneracy of our sesquilinear form, and the basic
identity\[r\varphi=\varphi z_r \text{ \qquad $\forall \, r \in R^W$,
$\forall\, \varphi \in M$}\] where $z_r$ denotes the element of the
center of~$H$ that corresponds to~$r$.
\end{proof}

\subsection{Restriction of the Kazhdan-Lusztig involution to the center}
Now consider the affine Hecke algebra ${\mathcal H}$ associated to
$G$.  This is an algebra over the ring ${\mathbb Z}[v,v^{-1}]$
($v$ an indeterminate), generated by symbols $T_w$ ($w$ ranging
over the extended affine Weyl group for $G$), which satisfy the
usual braid and quadratic relations.  If $q = p^n$ denotes the
cardinality of the residue field of $F$, the map $v \mapsto
q^{1/2}$ determines a ring homomorphism ${\mathbb Z}[v,v^{-1}] \to
\mathbf C$.  There is a canonical isomorphism $H = {\mathcal H}
\otimes_{{\mathbb Z}[v,v^{-1}]} {\mathbf C}$ (see section 7.2).

The Kazhdan-Lusztig involution $h
\mapsto \bar{h}$ of ${\mathcal H}$ is determined by $v \mapsto
v^{-1}$ and $T_w \mapsto T^{-1}_{w^{-1}}$ (beware that this does
not descend to an involution of $H$).  There is also an
anti-involution on ${\mathcal H}$ given by $v \mapsto v$ and $T_w
\mapsto T_{w^{-1}}$.  On specializing $v \to q^{1/2}$, this does
descend to $H$ and gives precisely the anti-involution of $H$
denoted $\iota$ above; therefore we denote the anti-involution of
${\mathcal H}$ by the same symbol.

 For each dominant coweight $\mu
\in X_*(A)$, we let $z_\mu = \sum_{\lambda \in W\mu}
\Theta_\lambda$, where $\Theta_\lambda$ is the element of
${\mathcal H}$ defined by $\Theta_\lambda =  v^{\langle 2\rho,
-\lambda_1 + \lambda_2 \rangle } T_{t_{\lambda_1}}
T^{-1}_{t_{\lambda_2}}$, where $\lambda = \lambda_1 - \lambda_2$,
and $\lambda_i$ is dominant ($i = 1,2$); see Remark \ref{rem172}.
 A result
of Bernstein says that the elements $z_\mu$ form a ${\mathbb
Z}[v,v^{-1}]$-basis for the center of ${\mathcal H}$, as $\mu$
ranges over dominant coweights in $X_*(A)$ (see also section
\ref{v-center}).

 These
considerations yield a simple proof of Corollary 8.8 in
\cite{lusztig83}:
\begin{lemma} \label{z-bar}
$\overline{z_\mu} = z_\mu$.
\end{lemma}

\begin{proof}

First of all we can relate the two involutions by the
easily-checked formula $\iota(\Theta_{-\lambda}) =
\overline{\Theta}_{\lambda}$.  It follows that $\iota (z_{-w_0\mu})
= \overline{z_\mu}$, where $w_0$ is the longest element of $W$.  On
the other hand, the previous lemma says that, at least after $v$ is
specialized to $q^{1/2}$, the elements $\iota(z_{-w_0\mu})$ and
$z_{\mu}$ coincide as elements in $H$.  Since this is true for
every power of $q$ by the same token, we must have the equality
$\iota(z_{-w_0\mu}) = z_\mu$ in ${\mathcal H}$ as well, which
proves that $\overline{z_{\mu}} = z_\mu$.
\end{proof}

\section{Satake isomorphism \cite{Sat}}
\subsection{Definition of $H_K$ and $M_K$}
Let $e_K$ be the idempotent $1_K/\meas(K)$
in~$H$. Put $H_K:=C_c(K\backslash G/K)$, which we identify with the subring
$e_K H e_K$ of~$H$ (so that $1_K \mapsto e_K$). We also put
$M_K:=C_c(A_\mathcal O N\backslash G/K)$, which we identify with the
$H_K$-submodule $Me_K$ of~$M$. Then $M_K$ is an $(R,H_K)$-bimodule, with
$R$-module structure inherited from the one on~$M$. Concretely, the action
of the function $h \in H_K$ on $m \in M_K$ is given by $m\ast h$, where
$\ast$ denotes convolution using the Haar measure on $G$ giving $K$
measure~$1$.

\subsection{The Satake transform}
Since $M_K$ is free of rank $1$ as $R$-module (with basis element the
spherical vector
$1_{A_\mathcal O N K}$), we get a $\mathbf C$-algebra homomorphism $H_K \to
R$, denoted $h \mapsto h^\vee$ and called the Satake transform,
characterized by the property that
\begin{equation} \label{sat1}
m \ast h= h^\vee \cdot m
\end{equation}
 for all $h \in H_K$ and all $m \in M_K$.

Taking $m$ to be the spherical vector, we get the equation
\begin{equation} \label{sat2}
1_{A_\mathcal O N K} \ast h= h^\vee \cdot 1_{A_\mathcal O N K}.
\end{equation}
In fact $h^\vee$ lies in the subalgebra $R^W$, as one sees by applying the
normalized intertwiners (which fix the spherical vector) to equation
\eqref{sat2}. Thus the Satake transform actually maps $H_K$ into $R^W$.

Recall that $\pi^\nu$ (for $\nu$ ranging through $X_*$) form a $\mathbf
C$-basis for $R$. Evaluating both sides of equation \eqref{sat2} on the
element $\pi^\nu$ and using the usual $G=ANK$ integration formula (see
\cite{cartier79}), one sees that the coefficient of $\pi^\nu$ in $h^\vee$ is
equal to
\begin{equation*}
\delta_B(\pi^\nu)^{-1/2} \int_N h(n\pi^\nu) \, dn,
\end{equation*}
where the Haar measure $dn$ is normalized so that $N_\mathcal O$ has
measure~$1$.

\subsection{Satake transform is an isomorphism}
The elements $h_\mu:=1_{K\pi^\mu K}$, with $\mu$ a dominant coweight, form a
$\mathbf C$-basis for $H_K$. The elements $s_\nu:=\sum_{\lambda \in W\nu}
\pi^{\lambda}$, with $\nu$ a dominant coweight, form a $\mathbf C$-basis
for~$R^W$. The coefficients $c_{\mu\nu}$ of $h_\mu^\vee$ in the basis
$s_\nu$ are given by
\begin{equation}
c_{\mu\nu}=\delta_B(\pi^\nu)^{-1/2} \int_N 1_{K\pi^\mu K}(n\pi^\nu) dn.
\end{equation}
The real number $c_{\mu\nu}$ is non-negative
and is non-zero if and only if $K \pi^\mu K$ meets $N \pi^\nu$. It follows
from
\cite[4.4.4]{BT} that $c_{\mu\nu}$ is $0$ unless $\nu \le \mu$ (by which we
mean that $\mu-\nu$ is a non-negative integral linear combination of simple
coroots), and it is obvious that $c_{\mu\mu}$ is non-zero. Therefore a
standard upper-triangular argument shows that the Satake transform is an
isomorphism from $H_K$ to $R^W$. In particular $H_K$ is commutative.

We remark that in \cite[Th\'{e}orem\`{e} 5.3.17]{Mats},  \cite{rap},
 \cite{haines'}, it is shown that if $\nu \le \mu$ (both dominant), then
$c_{\mu\nu}$ is non-zero.

\subsection{Compatibility of two involutions}
Recall from \ref{invdefs} the involutions $\iota$, $\iota_A$ on $H$,
$R$ respectively. It is clear that  $\iota$ preserves the subring
$H_K$ and that $\iota_A$ preserves the subring $R^W$. One sees easily
(imitate the proof of Lemma \ref{iota}) that the Satake isomorphism is
compatible with these involutions, in the sense that
\begin{equation}
(\iota(h))^\vee=\iota_A(h^\vee).
\end{equation}

\subsection{Further discussion of the Satake transform}
Consider a quasicharacter $\chi:A/A_\mathcal O \to \mathbf
C^\times$. Then $\chi$ determines a $\mathbf C$-algebra
homomorphism $R \to \mathbf C$. Using this homomorphism to extend
scalars, we obtain an $H_K$-module $\mathbf C \otimes_R M_K$ which  can be
identified with the ($1$-dimensional)   space of $K$-fixed vectors in the
unramified principal series representation
$i^G_B(\chi^{-1})$. It is customary to
work with left $G$-modules (and hence left modules over Hecke algebras) rather
than right modules, and to switch back and forth between right and left one
uses $g \mapsto g^{-1}$ on $G$ (and hence the involution $\iota$ on $H_K$).
Bearing these remarks in mind, one sees that for any $h \in H_K$ and any
$K$-fixed vector $v \in i^G_B(\chi)$ there is an equality
\begin{equation}\label{hec.chi}
hv=h^\vee(\chi)v,
\end{equation}
where for $r \in R$ we write $r(\chi)$ for the image of $r$ under the
homomorphism $R \to \mathbf C$ determined by $\chi$. (We used that
$\iota_A(r)(\chi^{-1})=r(\chi)$.)

\subsection{Compatibility of the Satake and Bernstein
isomorphisms \cite{dat, haines, lusztig83}}

We now have canonical isomorphisms (Satake and Bernstein)
\begin{equation*}
H_K \simeq R^W \simeq Z(H),
\end{equation*}
where $Z(H)$ denotes the center of~$H$. Let $h \in H_K$, $r \in R^W$, $z \in
Z(H)$ be elements that correspond to each other under these isomorphisms. We
have
\begin{equation*}
mh=me_Kh= r  me_K=me_Kz
\end{equation*}
for all $m \in M$.
It follows that
\begin{equation}
h=e_Kz,
\end{equation}
which is the compatibility referred to in the heading of this section.

\section{Macdonald's formula \cite{C',Macd,Mats}}
\subsection{Preliminary remarks about unramified matrix coefficients}
The contragredient of the induced representation $i^G_B(\chi)$ is
$i^G_B(\chi^{-1})$. (Recall that $i^G_B(\chi^{-1})$ has as usual a left $G$-action given by right translations.)  Now choose $K$-fixed vectors $v \in i^G_B(\chi)$ and
$\tilde v \in i^G_B(\chi^{-1})$ such that $\langle v,\tilde v \rangle =1$, and
put
\begin{equation}
\Gamma_\chi(g):= \langle gv,\tilde v \rangle ,
\end{equation}
an unramified matrix coefficient, otherwise known as a zonal spherical
function. Clearly $\Gamma_\chi$ is a $\mathbf C$-valued function on
$K\backslash G/K$, and we have
\begin{equation} \label{sph1}
\Gamma_\chi(1)=1.
\end{equation}
Let $h \in H_K$. It follows from the definition of $\Gamma_\chi$ that $(h \ast
\Gamma_\chi)(g)=\langle gv,h\tilde v \rangle$, which by \eqref{hec.chi} is
equal to
$h^\vee(\chi^{-1})\Gamma_\chi(g)$. Thus we have
\begin{equation}\label{sph2}
h \ast \Gamma_\chi = h^\vee(\chi^{-1}) \Gamma_\chi.
\end{equation}
Similarly we have
\begin{equation} \label{sph3}
 \Gamma_\chi \ast h = h^\vee(\chi^{-1}) \Gamma_\chi.
\end{equation}
The function $\Gamma_\chi$ is uniquely determined by \eqref{sph1} and either
of \eqref{sph2}, \eqref{sph3}; indeed, taking
$h=1_{K\pi^{-\mu}K}=\iota(K\pi^\mu K)$ in \eqref{sph2} and then evaluating
both sides  at the identity element, we see that
\begin{equation}\label{zsf=sat}
\meas(K\pi^\mu K)\cdot \Gamma_\chi(\pi^\mu)=(1_{K\pi^\mu K})^\vee(\chi),
\end{equation}
where the measure is taken with respect to the Haar measure on $G$ that gives
$K$ measure~$1$.
 In other words, knowing the values of unramified matrix coefficients is
essentially the same as knowing the Satake transforms of the elements $1_{K
\pi^\mu K}
\in H_K$.

\subsection{Definition of $\Gamma$}
It is more convenient to work with the $R$-valued matrix coefficient $\Gamma$
defined by
\begin{equation}\label{521}
\Gamma(g):=(1_{A_O NK}, 1_{A_O NK} \cdot g),
\end{equation}
where $(\cdot,\cdot)$ is our sesquilinear form on $i^G_B(\chi^{-1}_{\univ})$
(regarded as a right $G$-module).

 Of course $\Gamma$ is a function on $K\backslash
G/K$ with values in~$R$;  applying the homomorphism
$R \to \mathbf C$ determined by $\chi$ to the values of~$\Gamma$, we
get the $\mathbf C$-valued function $\Gamma_\chi$.
Therefore computing $\Gamma$ is the same as computing
$\Gamma_\chi$ for all $\chi$.

We can rewrite \eqref{521} as
\begin{equation}
\Gamma(g):=(1_{A_O NK}, 1_{A_O NK} \ast e_{KgK}),
\end{equation}
where $e_{KgK}$ denotes $\meas(KgK)^{-1}\cdot 1_{KgK}$, from which it follows
that
\begin{equation}\label{zsfw}
\Gamma(g)=(e_{KgK})^\vee,
\end{equation}
 in agreement with \eqref{zsf=sat}.
Equation \eqref{zsfw} shows that $\Gamma$ actually takes values in~$R^W$ and
hence that $\Gamma_{w\chi}=\Gamma_\chi$ for all $w \in W$.

Macdonald's formula \cite{C',Macd,Mats}
 is an explicit formula for $\Gamma_\chi$, which we will
now derive, following Casselman's method  \cite{C'}. As mentioned above, it
is the same to give an explicit formula for~$\Gamma$, and this is what we will
do.

\subsection{Decomposition of the spherical vector as a sum of eigenvectors}
As a first step
towards Macdonald's formula, we are going to decompose the spherical vector
$1_{A_\mathcal O NK} \in M$ as a sum of eigenvectors for the action of the
commutative subalgebra $R$ of $H$. This can only be done in $M_{\gen}$.

Recall that $v_1$ denotes the element $1_{A_\mathcal O NI} \in M$. The vector
$v_1$ is an eigenvector for the subalgebra $R$ of $H$ by the very definition
of that subalgebra; more precisely we have the formula
\begin{equation}\label{531}
v_1 \Theta_\lambda = \pi^\lambda \cdot v_1,
\end{equation}
where (as before) $\Theta_\lambda$ is a notation for the image of $\pi^\lambda
\in R$ under $R \hookrightarrow H$.
Applying the normalized intertwiner $K_w$ to this equation, we see that
\begin{equation}\label{532aaa}
K_w(v_1) \Theta_\lambda = \pi^{w\lambda} \cdot K_w(v_1),
\end{equation}
which shows that $K_w(v_1)$ is an eigenvector for~$R$ with character
$w^{-1}(\chi_{\univ})$.

\begin{lemma}\label{efnexp}
In $M_{\gen}$ we have the formula
\begin{equation*}
1_{A_\mathcal O NK}=\sum_{w \in W} w \Bigl( \prod_{\alpha > 0} \frac{1-q
\pi^{\alpha^\vee}}{1-\pi^{\alpha^\vee}} \Bigr) \cdot K_{w}(v_1).
\end{equation*}
\end{lemma}

\begin{proof} Let  $w_0$ denote the longest element of~$W$.
Recall the standard basis elements $v_x$ for $M$. Then $v_w$ ($w \in W$) form
an $R$-basis for $M$, hence an $L$-basis for $M_{\gen}$. From \eqref{triang}
it is clear that the vectors $K_w(v_1)$ also form an $L$-basis for $M_{\gen}$.
Write the spherical vector in this second basis:
\begin{equation}\label{dw}
1_{A_\mathcal O NK}=\sum_{w \in W} d_w \cdot  K_w(v_1).
\end{equation}
We can also write the spherical vector in the first basis; since the spherical
vector is equal to $\sum_{w \in W} v_w$, it is clear that the coefficient of
the basis element $v_{w_0}$ in the spherical vector is~$1$; on the other hand,
from \eqref{triang} and \eqref{dw}, it is clear that this same coefficient is
also equal to $d_{w_0}a_{w_0w_0}$; equating the two expressions for the
coefficient and using the explicit formula \eqref{diag.elts} for $a_{w_0w_0}$,
we see that
\begin{equation*}
d_{w_0}=\prod_{\alpha > 0}
 \frac{1-q\pi^{-\alpha^\vee}}{1-\pi^{-\alpha^\vee}}.
\end{equation*}
Moreover, since the normalized intertwiners $K_w$ fix the spherical vector, we
have
\begin{equation*}
d_{w_1w_2}=w_1(d_{w_2})
\end{equation*}
for all $w_1, w_2 \in W$, from which it follows that
\begin{equation*}
d_{w}=w \Bigl( \prod_{\alpha > 0} \frac{1-q
\pi^{\alpha^\vee}}{1-\pi^{\alpha^\vee}} \Bigr).
\end{equation*}
This completes the proof.
\end{proof}

\subsection{Partial information about some more matrix coefficients}
We see from Lemma \ref{efnexp} that in order to calculate $\Gamma$ it would be
enough to calculate the matrix coefficients $(1_{A_\mathcal O NK},K_w(v_1) \cdot
g)$. Now this new matrix coefficient is a function on $I\backslash G/K=X_*$
rather than $K\backslash G/K$, and it is difficult to calculate all its values.
Fortunately it easy to calculate them for elements $g$ of the form $\pi^\mu$ for
dominant coweights $\mu$, and in the end this is enough since $\Gamma$ is
$K$-bi-invariant and hence determined by its values on such elements.

\begin{lemma}\label{asym}
For group elements of the form $g=\pi^\mu$ with $\mu$ dominant, we have
\begin{equation*}
(1_{A_\mathcal O NK},K_w(v_1) \cdot g)=[K:I]^{-1}\cdot
\delta_B(\pi^\mu)^{1/2}
\cdot \pi^{w\mu}.
\end{equation*}
\end{lemma}

\begin{proof}
For $f \in M$, the function $f \cdot g$ need not be right $I$-invariant, 
and so need not have a simple form.  However, letting $\delta_g$ denote the
Dirac measure concentrated at $g$, we have an equality of measures $e_I
\cdot \delta_g \cdot e_I = e_{IgI}$, where $e_{X}$ is the characteristic
function of a set $X$ divided by its measure.  Since $g=\pi^\mu$ with $\mu$
dominant, we have
$e_{IgI}=\delta_B(\pi^\mu)T_{\pi^\mu}
=\delta_B(\pi^\mu)^{1/2}\Theta_\mu$.
Using these considerations (and the fact that the idempotent $e_I$ fixes
both $1_{A_\mathcal O NK}$ and $K_w(v_1)$), we see that 
\begin{equation*}
(1_{A_\mathcal O NK},K_w(v_1) \cdot
g)=\delta_B(\pi^\mu)^{1/2}\cdot(1_{A_\mathcal O NK},K_w(v_1)
\Theta_\mu).
\end{equation*}
{}From \eqref{532aaa} we have
$K_w(v_1)\Theta_\mu=\pi^{w\mu} K_w(v_1)$, and therefore
\begin{equation*}
(1_{A_\mathcal O NK},K_w(v_1) \cdot
g)=\delta_B(\pi^\mu)^{1/2}\cdot \pi^{w\mu} \cdot (1_{A_\mathcal O
NK},K_w(v_1)).
\end{equation*}
Using \eqref{223abc}, we see that
\begin{equation*}
(1_{A_\mathcal O NK},K_w(v_1))=w(K_{w^{-1}}(1_{A_\mathcal O NK}),v_1)=
w(1_{A_\mathcal O NK},v_1).
\end{equation*}
Moreover $(1_{A_\mathcal O NK},v_1)=[K:I]^{-1}$, as follows immediately from the
definitions. This completes the proof.

\end{proof}

\subsection{Macdonald's formula}
Combining Lemmas \ref{efnexp} and \ref{asym},  we
get
\begin{theorem}[Macdonald] \label{macdonald's}
For any  dominant coweight $\mu$ we have
\begin{equation*}
\Gamma(\pi^\mu) =[K:I]^{-1}\cdot\sum_{w \in W} w \Bigl( \prod_{\alpha > 0}
\frac{1-q\pi^{\alpha^\vee}}{1-\pi^{\alpha^\vee}} \Bigr) \cdot
\delta_B(\pi^\mu)^{1/2} \cdot \pi^{w\mu}
\end{equation*}
and
\begin{equation*}
\Gamma_\chi(\pi^\mu) =[K:I]^{-1}\cdot\sum_{w \in W} \Bigl( \prod_{\alpha > 0}
\frac{1-q
(w\chi)(\pi^{\alpha^\vee})}{1-(w\chi)(\pi^{\alpha^\vee})} \Bigr) \cdot
\delta_B(\pi^\mu)^{1/2} \cdot (w\chi)(\pi^{\mu}).
\end{equation*}
\end{theorem}

\subsection{Alternative version of Macdonald's formula}

For any finite subset $X \subset \widetilde W$, define the 
polynomial $X(t) := \sum_{w \in X}t^{l(w)}$, where the length function
$l(\cdot)$ is defined using the set of reflections for the
$\bar{B}$-positive simple affine roots, as in section 7.1.  Let $W_\mu$
denote the stabilizer of $\mu$ in $W$. We write $t_\mu$ for the element
$\pi^\mu$ of the translation subgroup of~$\widetilde W$. 

\begin{theorem}\label{alt.macdonald} For any dominant coweight $\mu$,
\begin{equation*}
(1_{K\pi^\mu K})^
\vee  = \frac{q^{\langle \rho, \mu \rangle}}{W_{\mu}(q^{-1})} \sum_{w \in W}
w \Bigl( \prod_{\alpha > 0}
\frac{1-q^{-1}\pi^{-\alpha^\vee}}{1-\pi^{-\alpha^\vee}} \Bigr) 
\cdot \pi^{w\mu}.
\end{equation*}
\end{theorem}

\begin{proof}
Using Theorem \ref{macdonald's}, this follows easily from the 
identities $[K:I] = W(q)$, $W(q) = q^{l(w_0)}W(q^{-1})$, ${\rm
meas}(K\pi^\mu K) = Wt_\mu W(q)/W(q)$, and
\begin{equation}\label{561aaa}
Wt_\mu W(q) = \frac{W(q)q^{l(t_\mu)}W(q^{-1})}{W_\mu(q^{-1})}.
\end{equation}
To prove (\ref{561aaa}), note that any element in $Wt_\mu W$ has 
a unique decomposition of the form $w^\mu  t_\mu w$, where $w \in W$ and
$w^\mu$ is a minimal length representative for a coset in $W/W_\mu$. 
Furthermore, such an element has length $l(w) + l(t_\mu) - l(w^\mu)$ (as may
be seen by induction on $l(w^\mu)$; note that $\mu$ is anti-dominant for
$\bar{B}$).
\end{proof}

\section{Casselman-Shalika formula \cite{CS,reeder',shintani}}

We are going to give an exposition of Casselman-Shalika's proof of their
formula \cite{CS} for unramified Whittaker functions.
\subsection{Unramified  characters $\psi$ on $\bar N$} Let $\Delta$ 
denote the set of simple roots. The abelian group
$\prod_{\alpha \in \Delta} N_{-\alpha}$ is a quotient of~$\bar N$. Here
$N_{-\alpha}$ is the root subgroup for $-\alpha$, which we identify with the
additive group $\mathbf G_a$ over $\mathcal O$.  Given characters
$\psi_\alpha:
 N_{-\alpha}
\to \mathbf C^\times$, their product defines a character on $\prod_{\alpha \in
\Delta} N_{-\alpha}$ and hence a character $\psi$ on~$\bar N$. We say that
$\psi$ is \emph{principal} if all the $\psi_\alpha$ are non-trivial. We say that
$\psi$ is \emph{unramified} if  all the characters
$\psi_\alpha$ are trivial on $\mathcal O$ but non-trivial on $\mathfrak
p^{-1}$.

\subsection{Whittaker functionals} \label{wh.func}
 Let $\psi$ be a principal character on~$\bar
N$.  Let
$S$ be a commutative
$\mathbf C$-algebra. The inclusion of $\mathbf C$ in $S$ lets us view $\psi$ as a
character with values in~$S^\times$.

Let $\chi:A \to S^\times$ be an $S$-valued character, and form the induced
representation $i^G_B(\chi)$, which is both a $G$-module and an $S$-module. A
\emph{Whittaker functional} on $i^G_B(\chi)$ is an $S$-module map
\begin{equation*}
L:i^G_B(\chi) \to S
\end{equation*}
such that $L(\bar n \phi)=\psi(\bar n)L(\phi)$ for all $\bar n \in \bar N$ and
all $\phi \in i^G_B(\chi)$.

In case $S=\mathbf C$ Rodier \cite{rodier} (see also \cite{CS}) proved that the
space of Whittaker functionals is $1$-dimensional and that there exists a unique
Whittaker functional $W$ whose restriction to the subspace of functions $\phi$ in
$i^G_B(\chi)$ supported on the big cell $B\bar N$ is given by the
integral
\begin{equation}\label{wh.int}
W(\phi)=\int_{\bar N} \phi(\bar n)\psi(\bar n)^{-1} \, d\bar n
\end{equation}
(the integrand of which is compactly supported by our assumption on the support
of $\phi$). Here $d\bar n$ denotes the Haar measure on $\bar N$ that gives
measure $1$ to $\bar N \cap K$.  For general
$S$ the same proof shows that there again exists a unique Whittaker functional
$W$ given by
\eqref{wh.int} for functions supported on the big cell, and that the $S$-module
of all Whittaker functionals is free of rank $1$ with $W$ as basis element.

Now let us consider the case in which $S$ is $R$ and $\chi$ is
$\chi_{\univ}^{-1}$. We let $J$ denote the set of negative coroots and consider the
completion $R_J$ of~$R$ defined in \ref{gener}. It follows from Lemma
\ref{compact} that the integral \eqref{wh.int} makes sense as an element of
$R_J$ for all $\phi \in i^G_B(\chi_{\univ}^{-1})$, and even for $\phi \in
i^G_B(\chi_{\univ}^{-1})\otimes_R R_J$. Using the uniqueness of $W$ for $R_J$,
we see that for $\phi \in i^G_B(\chi_{\univ}^{-1})$ the integral \eqref{wh.int}
actually takes values in the subring $R$ of $R_J$. (In other words the presence
of the principal character $\psi$ causes all but finitely many of the
coefficients of the Laurent power series $W(\phi)$ to vanish.) Therefore we will
now regard the Whittaker functional $W$ on $i^G_B(\chi_{\univ}^{-1})$
as being defined by the integral \eqref{wh.int}.

Recall that we have identified the module $M$ with the Iwahori-fixed vectors in
$i^G_B(\chi_{\univ}^{-1})$, and thus we also have the (restricted) Whittaker
functional $W:M \to R$. It is necessary to calculate $W$ for a few very special
vectors in~$M$.

From now on, we assume the character $\psi$ is principal and unramified.

\begin{lemma} \label{W.calc}
 Let $w_0$ denote the longest element in $W$, and let $\alpha$ be a
simple root with corresponding simple reflection $s_\alpha$. Then
\begin{gather}
W(v_1)=q^{-l(w_0)} \tag{i} \\
W(v_1+v_{s_\alpha})=q^{1-l(w_0)}\cdot (1-q^{-1}\pi^{-\alpha^\vee}).  \tag{ii}
\end{gather}
\end{lemma}
\begin{proof}
The first statement follows from the fact that $\bar N \cap BI =\bar N \cap I$,
which has measure $q^{-l(w_0)}$.  Similarly, the second statement
reduces to a calculation in
$SL(2)$, which we leave to the reader. Note that for $SL(2)$ the second statement
gives the value  (namely
$1-q^{-1}\pi^{-\alpha^\vee}$) of the Whittaker functional on the spherical
vector.

\end{proof}

\subsection{Effect of intertwiners on the Whittaker functional}
Earlier we defined normalized intertwiners $K_w$, normalized in the sense that
they preserve the spherical vector $1_{A_\mathcal O NK}$. Now we normalize them
differently. Put
\begin{equation}\label{K'}
K'_w:=\Bigl(\prod_{\alpha \in
R_w}\frac{1-q^{-1}\pi^{\alpha^\vee}}{1-q^{-1}\pi^{-\alpha^\vee}}\Bigr) \cdot
K_w=\Bigl(\prod_{\alpha \in
R_w}\frac{1-\pi^{\alpha^\vee}}{1-q^{-1}\pi^{-\alpha^\vee}}\Bigr)
\cdot I_w.
\end{equation}

\begin{lemma}[\cite{jacquet},\cite{CS}]\label{nor.wh}
The newly normalized intertwiners $K'_w$ preserve the
Whittaker functional $W$ in the sense that $W \circ K'_w=w \circ W$ for all $w
\in W$. On the right side of this equality
$w$ stands for the
automorphism of~$R$ determined by~$w$. Moreover $K'_{w_1w_2}=K'_{w_1}K'_{w_2}$.
\end{lemma}
\begin{proof}
One sees directly from the definition that $K'_w$ is multiplicative in~$w$.
Therefore to prove the first statement of the lemma, it is enough to treat the
case $w=s_\alpha$ for a simple root $\alpha$. By uniqueness of~$W$ there exists
$c \in L^\times$ such that
\begin{equation}\label{find.c}
W \circ K'_{s_\alpha}=c (s_{\alpha}
\circ W). \footnote{Here, we are implicitly using normalized intertwiners 
$K_w : L \otimes_R i^G_B(\chi^{-1}_{\univ}) \rightarrow L \otimes_R i^G_B(\chi^{-1}_{\univ})$.   The reader may derive the existence and basic properties of 
such 
intertwiners following 
the method of sections 1.10-2.2. 
Although the discussion there was limited to the theory of 
intertwiners on the Iwahori-invariants 
in the induced modules in question, it is  
possible to develop a similar theory on 
the induced modules themselves.}
\end{equation}
To prove that $c=1$ we evaluate both sides of \eqref{find.c} on
$v_1+v_{s_\alpha}$, using  Lemma \ref{lemma4.1}(ii) and Lemma \ref{W.calc}(ii).
\end{proof}

\begin{lemma}\label{efnexp'}
In $M_{\gen}$ we have the formula
\begin{equation*}
1_{A_\mathcal O NK}=q^{l(w_0)}\cdot \Bigl(\prod_{\alpha >
0}(1-q^{-1}\pi^{-\alpha^\vee})
\Bigr)
\cdot
\sum_{w
\in W} w
\Bigl(
\prod_{\alpha > 0}
\frac{1}{1-\pi^{-\alpha^\vee}} \Bigr) \cdot K'_{w}(v_1).
\end{equation*}
\end{lemma}
\begin{proof}
This follows from Lemma \ref{efnexp} and the definition of~$K'_w$.
\end{proof}

\subsection {Whittaker functions} We continue with $S$, $\psi$ and $\chi$ as in
\ref{wh.func}. For any $\phi \in i^G_B(\chi)$ we define the
corresponding Whittaker function $\mathcal W_\phi: G \to S$ by
\begin{equation}
\mathcal W_\phi(g):=W(g\phi).
\end{equation}
Then $\mathcal W_\phi$ is an $S$-valued function satisfying the transformation
law
\begin{equation} \label{tr.law}
f(\bar ng)=\psi(\bar n) f(g)  \qquad \forall \,\, \bar n \in \bar N,
\end{equation}
and $\phi \mapsto \mathcal W_\phi$ is a $G$-map from $i^G_B(\chi)$ to the space
of functions satisfying \eqref{tr.law}.

\subsection {Unramified Whittaker functions} From now on we assume that both
$\phi$ and $\chi$ are unramified. Then inside $i^G_B(\chi)$ we have the
normalized spherical vector $\phi_{\chi}$ defined by
\begin{equation*}
\phi_{\chi}(ank)=\chi(a)\delta_B(a)^{1/2}.
\end{equation*}
The Casselman-Shalika formula is an explicit formula for the Whittaker function
$\mathcal W_\chi:=\mathcal W_{\phi_{\chi}}$ corresponding to the spherical vector
$\phi_{\chi}$. It is enough to consider the case in which $S$ is $R$ and
$\chi$ is $\chi^{-1}_{\univ}$, in which case we abbreviate $\mathcal
W_{\chi_{\univ}^{-1}}$ to $\mathcal W$.

Since $\mathcal W$ is right $K$-invariant and satisfies \eqref{tr.law}, it is
determined by its values on elements $g \in G$ of the form $\pi^{-\mu}$ for $\mu
\in X_*$. In fact $\mathcal W(\pi^{-\mu})=0$ unless $\mu$ is dominant. Indeed,
for
$x \in N_{-\alpha} \cap K$ we have
\begin{equation*}
\mathcal W(\pi^{-\mu})=\mathcal
W(\pi^{-\mu}x)=\psi_\alpha(\pi^{-\mu}x\pi^\mu)\mathcal W(\pi^{-\mu}),
\end{equation*}
which implies that $\mathcal W(\pi^{-\mu})$ vanishes unless $\psi_\alpha$ is
trivial on $\mathfrak p^{\langle \alpha,\mu \rangle}$, which, since $\psi$ is
unramified, implies in turn that $  \langle \alpha,\mu \rangle \ge 0$. Therefore
it is enough to find the values of $\mathcal W(g)$ for $g$ of the form
$\pi^{-\mu}$ for dominant $\mu$.

\begin{theorem}[Casselman-Shalika]
Let $\mu$ be a dominant coweight. Then
\begin{equation*}
\mathcal W(\pi^{-\mu})=\bigl ( \prod_{\alpha > 0} (1 -
q^{-1}\pi^{-\alpha^\vee})\bigr )
\cdot \delta_B(\pi^\mu)^{1/2} \cdot E_\mu,
\end{equation*}
where $E_\mu \in R^W$ is the character of the irreducible representation of the
Langlands dual group $G^\vee$ having highest weight $\mu$.
\end{theorem}

\begin{proof}
As usual (see \ref{sec.pt.view}) we identify  $i^G_B(\chi^{-1}_{\univ})$
with  $C^\infty_c(A_{\mathcal O} N\backslash G)$. The spherical vector
in the induced representation corresponds to
$1_{A_\mathcal O NK}$.  We begin by noting that
\begin{equation}
\mathcal W(\pi^{-\mu})=W(1_{A_\mathcal O NK\pi^\mu})=
W(1_{A_\mathcal O NK}\cdot
e_{I\pi^\mu I})
\end{equation}
where $e_{I\pi^\mu I}$ denotes the characteristic function of $I\pi^\mu I$
divided by its measure. Here we used 
that $I\pi^\mu I=I\pi^\mu (I\cap \bar N)$ (a
consequence of the Iwahori factorization $I=(I \cap B) 
\cdot (I\cap \bar N)$ and
the dominance of $\mu$), as well as the right $I$-invariance 
of $1_{A_{\mathcal
O}NK}$ and the fact that $\psi$ is trivial on $I \cap \bar N$.  Since
$e_{I\pi^\mu I}=\delta_B(\pi^\mu)^{1/2} \cdot \Theta_\mu$,
 we can rewrite the
equation above as
\begin{equation}
\mathcal W(\pi^{-\mu})=\delta_B(\pi^\mu)^{1/2} W(1_{A_\mathcal O NK}\cdot
\Theta_\mu).
\end{equation}
It then follows from Lemmas \ref{W.calc}(i), \ref{nor.wh}, and \ref{efnexp'}
that
\begin{equation*}
\mathcal W(\pi^{-\mu})= \delta_B(\pi^\mu)^{1/2}  \cdot \Bigl(\prod_{\alpha >
0}(1-q^{-1}\pi^{-\alpha^\vee})
\Bigr)
\cdot
\sum_{w
\in W} w
\Bigl(
\prod_{\alpha > 0}
\frac{1}{1-\pi^{-\alpha^\vee}} \Bigr) \cdot \pi^{w\mu}.
\end{equation*}

The Casselman-Shalika formula now follows from the Weyl character formula:

\begin{equation*}
E_\mu= \sum_{w
\in W} w
\Bigl(
\prod_{\alpha > 0}
\frac{\pi^\mu}{1-\pi^{-\alpha^\vee}} \Bigr) .
\end{equation*}
\end{proof}

\section{The Lusztig-Kato formula  \cite{kato, lusztig83}}\label{Kato}

Following the strategy of \cite{kato}, we derive the formula of Lusztig-Kato
and, as a corollary, another result of Lusztig \cite{lusztig83}.  The
Lusztig-Kato formula relates the Satake transforms of the functions
$1_{K\pi^\lambda K}$ to the character $E_\mu$ of the highest weight module
of the Langlands dual group corresponding to $\mu$.  It is the
function-theoretic counterpart of the geometric Satake isomorphism
\cite{Ginzburg, MV}, and it can also be formally deduced from that statement
by using the function-sheaf dictionary.

The proof requires us to give $v$-analogs of several objects studied above
(here $v$ is an indeterminate which can be specialized to $q^{1/2}$).  Most
importantly, we need the $v$-analog of Theorem \ref{alt.macdonald}.  The
reader willing to accept that on faith may skip directly to section 7.7.

\subsection{Preliminaries about affine roots} Write $T$ 
for the group $X_*(A)$, viewed as the group of translations in the
extended affine Weyl group $\widetilde{W}$; thus $\widetilde{W} = T \rtimes
W$. We denote by $t_\mu$ the element of $T$ corresponding to the
cocharacter $\mu$. For simplicity, we assume here that the root system
underlying
$G$ is irreducible. Let
$\alpha_1,
\dots ,
\alpha_r$ denote the $B$-positive simple roots, and let $\tilde \alpha$
denote the $B$-highest root.  Let $s_0 = t_{-\tilde{\alpha}^\vee}s_{\tilde
\alpha}$, and  $S_{\rm aff} = S \cup \{ s_0 \}$.  Here $S = \{s_{\alpha_i} =
s_{-\alpha_i}\}_{i=1}^r$ is the set of simple reflections corresponding to
the $B$-positive (or $\bar{B}$-positive) simple roots, but our definition of
$s_0$ means that $S_{\rm aff}$ is the set of simple affine reflections
corresponding to the $\bar{B}$-positive affine roots.

  We have $\widetilde{W} = W_{\rm aff} \rtimes \Omega$, where $W_{\rm
aff}$ is the Coxeter group generated by $S_{\rm aff}$, and $\Omega$ is the
subgroup of $\widetilde{W}$ which preserves the set of
$\bar{B}$-positive  simple affine roots under the usual left action (an
affine-linear automorphism acts on a functional by pre-composition with its
inverse).  The set $S_{\rm aff}$ induces a length function and a Bruhat
order on $\widetilde{W}$ (the same as that mentioned in Lemma \ref{rank1}).  The
elements $\sigma \in \Omega$ are of length zero, and the algebra generated
by the functions $1_{I\sigma I}$ is naturally isomorphic to ${\mathbf
C}[\Omega]$.  We have a twisted tensor product decomposition $H = H_{\rm
aff} \otimes {\mathbf C}[\Omega]$, where $H_{\rm aff}$ is the algebra
generated by the functions $1_{IxI}$, $x \in W_{\rm aff}$ (this follows from
the remarks following Lemma \ref{H-action} below).

Recall our convention for embedding $X_*(A)$ into $A$: $\lambda \mapsto
\pi^\lambda = \lambda(\pi)$.  We also regard each $w \in W$ as an element in
$K$, fixed once and for all.  These conventions tell us how we view elements
of $\widetilde W$ as elements in $G$.  For example, for ${\rm SL}(2)$, we are
identifying $s_0$ with the element $\begin{bmatrix} 0 & \pi^{-1} \\ -\pi & 0
\end{bmatrix}$.  It is important to bear these conventions in mind in this
section.

\subsection{More about the $H$-action on $M$}
\begin{lemma} \label{H-action} Fix $\varphi \in M$ and $\pi^\lambda w \in
\widetilde W$, where $w \in W$.  Suppose $\sigma \in \Omega$, and that $s = s_\alpha \in S$
corresponds to a $B$-positive simple root $\alpha$. Then we have
\begin{align*}
\varphi T_{s_\alpha}(\pi^\lambda w) &= \begin{cases} q \cdot
\varphi(\pi^\lambda w s), \,\, ~ \mbox{ if $w(\alpha)$ is $B$-positive} \\
\varphi(\pi^\lambda ws) + (q-1) \cdot \varphi(\pi^\lambda w), \,\, ~
\mbox{if $w(\alpha)$ is $B$-negative} \end{cases} \tag{i} \\
\varphi T_{s_0}(\pi^\lambda w) &= \begin{cases} q \cdot \varphi(\pi^\lambda
w s_0), \,\, ~ \mbox{ if $w(\tilde \alpha)$ is $B$-negative}. \\
\varphi(\pi^\lambda ws_0) + (q-1) \cdot \varphi(\pi^\lambda w), \,\, ~
\mbox{if $w(\tilde \alpha)$ is $B$-positive}. \end{cases} \tag{ii} \\
\varphi T_{\sigma}(\pi^\lambda w) &= \varphi(\pi^\lambda w \sigma^{-1}). 
\tag{iii}
\end{align*}
\end{lemma}

\begin{proof} To illustrate the method, we prove (ii).  Let $\{ x_i
\}_{i=0}^{q-1}$ denote a set of representatives for $\mathcal O / P$ taken
in $\mathcal O ^\times \cup \{0\}$.  For a root $\beta$, let $u_\beta:
{\mathbf G}_a \rightarrow G$ denote the associated homomorphism.  Then we
have the decomposition
\begin{equation*} Is_0 I = \coprod_{i} u_{-\tilde \alpha}(\pi x_i) s_0 I.
\end{equation*} We therefore have
\begin{equation*}
\varphi T_{s_0}(\pi^\lambda w) = \sum_i \varphi(\pi^\lambda w \, u_{-\tilde
\alpha}(\pi x_i) s_0).
\end{equation*} If $w(\tilde \alpha)$ is $B$-negative, then each term in the
sum is $\varphi(\pi^\lambda w s_0)$. If $w(\tilde \alpha)$ is $B$-positive,
then the term for $x_i = 0$ is $\varphi(\pi^\lambda w s_0)$.  If $x_i \neq
0$, then using the identity
\begin{equation*} u_{-\tilde \alpha}(\pi x_i) s_0 I = u_{\tilde
\alpha}(\pi^{-1} x^{-1}_i) I
\end{equation*} (which holds whenever $x_i \in \mathcal O ^\times$), we see
the term indexed by $x_i$ is
$\varphi(\pi^\lambda w)$. 

Part (i) can be proved in a similar way; alternatively it can be derived from
\eqref{use45} together with the usual relations in the Hecke algebra for
$W$. 
\end{proof}

The proof of (ii) above parallels the standard proof of the
Iwahori-Matsumoto relations in
$H$, which state that for
$x\in \widetilde W$, $s \in S_{\rm aff}$, and $\sigma \in \Omega$
\begin{align*} T_x T_s &= \begin{cases} T_{xs}, \,\,\, \mbox{if $x < xs$} \\
q \cdot T_{xs} + (q-1)\cdot T_x,
\,\,\, \mbox{if $xs < x$} \end{cases} \\ T_x T_\sigma &= T_{x\sigma},
\end{align*} where $<$ denotes the Bruhat order determined by $S_{\rm
aff}$.  If $\mathcal H$ denotes the affine Hecke algebra over ${\mathbb Z}_v
:= {\mathbb Z}[v,v^{-1}]$ associated to our root system, this means we have
a canonical isomorphism $H = \mathcal H
\otimes_{{\mathbb Z}_v} {\mathbf C}$.

Note that $R$ is the Iwahori-Hecke algebra for the group $A$, and hence it
also has a $v$-analog over ${\mathbb Z}_v$, which we denote by ${\mathcal
R}$.  Concretely, we have ${\mathcal R} = {\mathbb Z}_v[X_*(A)]$.

We will use Lemma \ref{H-action} as the starting point in defining
$v$-analogs
$$ {\mathcal M}, \hspace{.1in} i^{\widetilde W}_T(\chi^{-1}_{\univ}),
\hspace{.1in} ({\mathcal R}, {\mathcal H})- \, \mbox{actions}, 
\hspace{.1in} ( \cdot | \cdot ), \hspace{.1in} K_w, \hspace{.1in} \mathcal
M' e_W, \hspace{.1in} h^\vee
$$ of the objects we have already studied
$$ M, \hspace{.1in}  i^G_B(\chi^{-1}_{\univ})^I, \hspace{.1in} (R,H)-\,
\mbox{actions}, \hspace{.1in} ( \cdot , \cdot ), \hspace{.1in}  K_w,
\hspace{.1in}  M_K, \hspace{.1in} h^\vee.
$$

\subsection{$v$-analogs of $M$ and $i^G_B(\chi^{-1}_{\rm univ})^I$}

Let us define ${\mathcal M}$ to be the set of functions $\varphi: \widetilde W
\rightarrow {\mathbb Z}_v$ which are supported on a finite subset.  This is
a free ${\mathbb Z}_v$-module with basis given by the characteristic
functions $1_x$, $x \in \widetilde W$.

Next we define $\delta : T \rightarrow {\mathbb Z}_v^\times$ by
$\delta(t_\lambda):= v^{-2 \langle \rho, \lambda \rangle}$.  This is the
$v$-analog of the function $\delta_B$. By $\delta^{1/2}$ we will mean the
obvious square root of $\delta$, namely the character $t_\lambda \mapsto
v^{-\langle \rho,\lambda \rangle}$. 

The left action of ${\mathcal R}$ on ${\mathcal M}$ is given by the formula
$t \cdot 1_x : = \delta^{1/2}(t) 1_{tx}$.  The right ${\mathcal H}$-action
is given by defining (following Lemma \ref{H-action}), for $\varphi \in
{\mathcal M}$,
\begin{align*}
\varphi T_{s_\alpha}(t_\lambda w) &= \begin{cases} v^2 \cdot
\varphi(t_\lambda w s), \,\, ~ \mbox{ if $w(\alpha)$ is $B$-positive} \\
\varphi(t_\lambda ws) + (v^2-1) \cdot \varphi(t_\lambda w), \,\, ~ \mbox{if
$w(\alpha)$ is $B$-negative} \end{cases} \tag{i} \\
\varphi T_{s_0}(t_\lambda w) &= \begin{cases} v^2 \cdot \varphi(t_\lambda w
s_0), \,\, ~ \mbox{ if $w(\tilde \alpha)$ is $B$-negative}. \\
\varphi(t_\lambda ws_0) + (v^2-1) \cdot \varphi(t_\lambda w), \,\, ~
\mbox{if $w(\tilde \alpha)$ is $B$-positive}. \end{cases} \tag{ii} \\
\varphi T_{\sigma}(t_\lambda w) &= \varphi(t_\lambda w \sigma^{-1}). 
\tag{iii}
\end{align*} These rules determine a right ${\mathcal
H}$-module structure on ${\mathcal M}$ 
\footnote{Matsumoto gives a similar definition for a left ${\mathcal H}$-action in \cite[section 4.1.1]{Mats}.}, 
and moreover ${\mathcal M}$ is an $({\mathcal R}, {\mathcal H})$-bimodule. 
Indeed, it suffices to observe
that by Lemma
\ref{H-action} 
 this statement holds after every specialization $v \mapsto
q^{1/2}$. Specialization arguments like this will be
used repeatedly below to prove $v$-analogs 
of statements known for $M$.

Now define $i^{\widetilde W}_T(\chi^{-1}_{\univ})$ to be the set of functions
$\phi: \widetilde W \rightarrow {\mathcal R}$ which satisfy
\begin{equation*}
\phi(t x) = \delta^{1/2}(t) \phi(x) \cdot t^{-1},
\end{equation*} for $t \in T$ and $x \in \widetilde W$.  As in section 1.5,
there is a canonical isomorphism
\begin{equation} \label{721aaa} {\mathcal M} = i^{\tilde
W}_T(\chi^{-1}_{\univ}).
\end{equation} Explicitly, we associate to $\varphi \in {\mathcal M}$ the
function $\phi$ given by
$$
\phi(x) = \sum_{t \in T} \delta^{-1/2}(t) \,\, \varphi(tx) \cdot t.
$$ The left action of ${\mathcal R}$ on $i^{\widetilde W}_T(\chi^{-1}_{\univ})$
is defined by $(t\cdot \phi)(x) = t(\phi(x))$.  The right action of
${\mathcal H}$ is defined by requiring the isomorphism ${\mathcal M}
\rightarrow i^{\widetilde W}_T(\chi^{-1}_{\univ})$ to be ${\mathcal H}$-linear
(one could also write out an explicit rule, again in the spirit of Lemma
\ref{H-action}).  Then ${\mathcal M} = i^{\widetilde W}_T(\chi^{-1}_{\univ})$ is
an isomorphism of $({\mathcal R},{\mathcal H})$-bimodules.

\subsection{$v$-analog of the sesquilinear pairing}

We will define an ${\mathcal R}$-valued sesquilinear pairing $( \cdot |
\cdot )$ on $i^{\widetilde W}_T(\chi^{-1}_{\univ})$ (thus on ${\mathcal M}$)
which is {\em almost} the $v$-analog of $( \cdot , \cdot )$ (they differ by
a constant).  Hence, it will automatically satisfy the analogs of
(\ref{sesq}), (\ref{herm}), and (\ref{h-inv}). We write $\iota_\mathcal R$
for the $v$-analog of $\iota_R$, namely the  involution
on $\mathcal R=\mathbb Z_v[X_*(A)]$ induced by the identity on $\mathbb
Z_v$ and the map $\mu \mapsto -\mu $ on $X_*(A)$. 

For $\phi_1, \phi_2 \in i^{\widetilde W}_T(\chi^{-1}_{\univ})$, define
\begin{equation} \label{v-pairing} ( \phi_1 | \phi_2 ) = \sum_{w \in W}
v^{2l(w)} \,\, \iota_{\mathcal R} \phi_1(w) \, \phi_2(w).
\end{equation}

\begin{lemma} \label{almost_pairing}  The pairing $( \cdot | \cdot )$ on
${\mathcal M}$  induces the pairing $W(q)\, ( \cdot , \cdot )$ on $M =
{\mathcal M} \otimes_{{\mathbb Z}_v} {\mathbf C}$.
\end{lemma}

\begin{proof} Consider the ${\mathcal R}$-basis $\{ 1_w \}_{w \in W}$ for
${\mathcal M}$, and the corresponding  $R$-basis $\{ v_w \}_{w \in W}$ for 
$M$.   From the definitions, we easily see
\begin{equation*} ( 1_w \, | \, 1_{w'} ) = v^{2l(w)} \, \delta_{w,w'}.
\end{equation*} It is therefore enough to prove
\begin{equation*} ( v_w , v_{w'} ) = q^{l(w)} \, W(q)^{-1} \, \delta_{w,w'}.
\end{equation*} The orthogonality is clear, and then one can easily check
that
$$ (v_w , v_w) = (1_{A_{\mathcal O} N K} , v_w) = (1_{A_{\mathcal O} N K}
T_{w^{-1}} , v_1) = q^{l(w)} W(q)^{-1}.
$$
\end{proof}

\subsection{$v$-analogs of normalized intertwiners}\label{v-center}

For a simple reflection $s = s_\alpha$, define $J_{s} : {\mathcal M}
\rightarrow {\mathcal M}$ by setting
\begin{align*} J_s(1_1) &= v^{-2}(1 - t_{\alpha^\vee}) \cdot 1_s + (1 -
v^{-2})t_{\alpha^\vee} \cdot 1_1, \\ J_s(1_1 h) &= J_s(1_1)h, \,\,\,
\mbox{for $h \in {\mathcal H}$}.
\end{align*} This makes sense because ${\mathcal M}$ is the free ${\mathcal
H}$-module generated by $1_1$ (by the same upper-triangular argument we used
to prove that $M$ is the free $H$-module generated by $v_1$).   Further, for
any
$w
\in W$, choose a reduced expression $w = s_1 \cdots s_n$, and set
$$ J_w := J_{s_1}\circ \cdots \circ J_{s_n}.
$$ (The usual specialization argument shows that the right hand
side is independent of the choice of reduced expression.)  

Next define
${\mathcal L}$ to be the fraction field of ${\mathcal R}$; note that
$\mathcal L^W$ is the fraction field of~$\mathcal R^W$ and that $\mathcal
L=\mathcal L^W\otimes_{\mathcal R^W} \mathcal R$. Imitating what we did
before, we see that $\mathcal R^W$ embeds into the center of~$\mathcal H$,
so that we can form the algebra ${\mathcal H}_{\gen} := {\mathcal L^W}
\otimes_{\mathcal R^W} {\mathcal H}$ and the right ${\mathcal
H}_{\gen}$-module 
${\mathcal M}_{\gen} := {\mathcal L^W} \otimes_{\mathcal R^W} {\mathcal
M}={\mathcal L} \otimes_{\mathcal R} {\mathcal
M}$.  Finally we define the normalized intertwiner $K_w : {\mathcal
M}_{\gen}
\rightarrow {\mathcal M}_{\gen}$ by
\begin{equation*} K_w:=\Bigl(\prod_{\alpha \in
R_w}\frac{1}{1-v^{-2}t_{\alpha^\vee}}\Bigr) \cdot J_w.
\end{equation*}

It is clear that
\begin{enumerate}
\item[(i)] $K_w$ is $\,\, {\mathcal H}_{\gen}$-linear;
\item[(ii)] $K_w \circ t_\lambda = t_{w\lambda} \circ K_w$;
\item[(iii)] $K_w$ fixes $1_W$,
\end{enumerate} where $1_W := \sum_{w \in W} 1_w$ is the $v$-analog of
$1_{A_\mathcal O N K}$.   It is also clear that $w \mapsto K_w$ defines a
homomorphism $W \rightarrow {\mathcal H}_{\gen}^\times$, and that the
$v$-analogs of (2.2.3-2.2.5) hold (using $( \cdot | \cdot )$ in (2.2.3)).  
Moreover, the $v$-analog of Lemma \ref{efnexp} holds. Finally,  we recover
Bernstein's result that 
$\mathcal R^W$ is the center of $\mathcal H$ by going through all the
same steps we did before.

\subsection{$v$-analog of the Satake isomorphism}

Let ${\mathbb Z}'_v$, ${\mathcal R}'$, ${\mathcal H}'$, and ${\mathcal M}'$
denote the localizations of
${\mathbb Z}_v$, ${\mathcal R}$, ${\mathcal H}$, and ${\mathcal M}$ at the
element
$W(v^2) \in {\mathbb Z}_v$.  Let $T_W = \sum_{w \in W} T_w$ and $e_W =
W(v^2)^{-1}T_W$, an element in ${\mathcal H}'$.  Further, define ${\mathcal
H}_0 = e_W {\mathcal H}' e_W$, and
${\mathcal M}_0 = {\mathcal M}'e_W$.  Then ${\mathcal
H}_0$ is a
${\mathbb Z}'_v$-algebra with product $* := W(v^2)^{-1}\cdot$ and identity
element
$T_W$,  where $\cdot$ denotes the usual product in ${\mathcal H}$.  
Similarly ${\mathcal M}_0$ is an
$\mathcal H_0$-module with product $* := W(v^2)^{-1}\cdot$,  where $\cdot$
now denotes the 
usual ${\mathcal H}$-action on ${\mathcal M}$.  It is clear that $*$ makes
${\mathcal M}_0$ an
$({\mathcal R}',{\mathcal H}_0)$-bimodule.

The ${\mathcal R}'$-module ${\mathcal M}_0$ is free of rank 1, so there is a
homomorphism
$$
\vee : {\mathcal H}_0 \rightarrow {\mathcal R}'
$$ characterized by
\begin{equation*} m_0 * h_0 = h^{\vee}_0 m_0,
\end{equation*} for all $h_0 \in {\mathcal H}_0$ and all $m_0 \in {\mathcal
M}_0$.

We have the formula
\begin{equation}\label{761aaa} h^{\vee}_0 = W(v^2)^{-1}( 1_W \, | \, 1_W *
h_0 ) = ( 1_1 \, | \, 1_1h_0).
\end{equation}

We can now easily derive the $v$-analog of Theorem \ref{alt.macdonald}.  We
apply the first equality of  (\ref{761aaa}) to the function
\begin{equation} h_\mu := \sum_{w \in Wt_\mu W}T_w = 
\frac{W(v^2)W(v^{-2})}{W_\mu(v^{-2})} \,\, e_W T_{t_\mu}e_W
\end{equation} to get
\begin{align*} (h_\mu)^\vee &=  \frac{W(v^{-2})}{W(v^2)W_\mu(v^{-2})} \,\, (
1_W \, | \, 1_W \cdot e_W T_{t_{\mu}}e_W) \\ &= \frac{v^{-
2l(w_0)}v^{2(l(t_\mu)/2)}}{W_\mu(v^{-2})} \,\, (1_W \, | \, 1_W \Theta_\mu).
\end{align*} Now using the $v$-analog of Lemma \ref{efnexp} as in the proof
of Theorem \ref{macdonald's}, we find
\vspace{.2in}
\begin{equation} \label{v-alt.macdonald} (h_\mu)^\vee =
\frac{v^{2(l(t_\mu)/2)}}{W_\mu(v^{-2})}
\sum_{w \in W} w \Bigl( \prod_{\alpha > 0}
\frac{1-v^{-2}t_{-\alpha^\vee}}{1-t_{-\alpha^\vee}} \Bigr) \cdot t_{w\mu}.
\end{equation}

\subsection{The Satake isomorphism commutes with the Kazhdan-Lusztig
involution}

The compatibility of the Bernstein and Satake isomorphisms (4.6) is the
commutativity of the following diagram:
$$
\xymatrix{ {\mathcal R}'^W \ar[d]_B & e_W{\mathcal H}' e_W \ar[l]_{b} \\
Z({\mathcal H}') \ar[ur]_{- \cdot e_W},}
$$ where $b(h_0) := W(v^2)h^\vee_0$ and where the Bernstein isomorphism $B$
sends
$\sum_{\lambda \in W\mu}t_{\lambda}$ to $z_\mu$. By Lemma
\ref{z-bar}, $B$ commutes with the Kazhdan-Lusztig involution. Since
$\overline{e_W} = e_W$ 
\footnote{Via the function-sheaf dictionary, the Kazhdan-Lusztig involution corresponds to taking the Verdier dual.  The equality can then be derived from the fact that if the constant sheaf on the smooth variety $G/B$ is placed in degree  $-l(w_0)$ and Tate-twisted by $l(w_0)/2$, the resulting complex is Verdier self-dual.}
, the diagonal map does as well.  We thus have the
following lemma which is implicit in \cite{lusztig83}, section 8.

\begin{lemma} \label{KLcommuteswithSat} For every $h_0 \in {\mathcal H}_0$,
$$ b(\overline{h_0}) = \overline{b(h_0)}.
$$ Equivalently,
$$ (\overline{h_0})^\vee = v^{-2 \, l(w_0)} \, \overline{h^\vee_0}.
$$
\end{lemma}

Note that the Kazhdan-Lusztig involution on the commutative ring ${\mathcal
R}'$ is simply the map sending
$\sum_{\lambda}z_\lambda(v,v^{-1}) t_\lambda$ to $\sum_{\lambda}
z_\lambda(v^{-1},v) t_\lambda$.

\subsection{The Lusztig-Kato formula}

We now switch notation and let $q^{1/2}$ play the role of the indeterminate
$v$ used in sections 7.1-7.7. In this section we will use some elementary
properties of the Kazhdan-Lusztig polynomials $P_{x,y}(q)$ attached to $x,y
\in \widetilde W$, all of which may be found in \cite{KL}.

Recall that throughout this article, the Bruhat order $\leq$ and the length
function $l(\cdot)$ on $\widetilde W$ are defined using the $\bar{B}$-positive
affine reflections $S_{\rm aff} := S \cup \{s_0 \}$. For any dominant
coweight $\lambda$, the element $w_\lambda:= t_\lambda w_0$ is the unique
longest element in
$Wt_\lambda W$, and $l(t_\lambda w_0) = l(w_0) + l(t_\lambda) = l(w_0) +
2\langle \rho, \lambda \rangle$.  It is known that
\begin{equation*}
\{ x \leq w_\mu \} = \cup_{\lambda \preceq \mu} W t_\lambda W,
\end{equation*} where $\lambda$ ranges over dominant coweights such that
$\mu - \lambda$ is a sum of $B$-positive coroots.

\begin{theorem}[Lusztig, Kato] \label{Lusztig-Kato} For any dominant
coweight $\mu$, let $E_\mu$ denote the character of the corresponding
highest weight module of the Langlands dual group $G^\vee$.  Let $h_\mu$
denote the function $\sum_{w \, \in \, W t_\mu W} T_w$.  Then we have
\begin{equation*} E_{\mu} = \sum_{\lambda \preceq \mu} q^{-l(t_\mu)/2}
P_{w_\lambda, w_\mu}(q) \,\, (h_\lambda)^\vee.
\end{equation*}
\end{theorem}

\begin{proof} We have the identity
\begin{equation*}
\overline{q^{-l(y)/2} \sum_{x \leq y}P_{x,y}(q) T_x} = q^{-l(y)/2}\sum_{x
\leq y}P_{x,y}(q)T_x.
\end{equation*} Applying this to $y = w_\mu$ and using $P_{w w_\lambda w',
w_\mu}(q) = P_{w_\lambda, w_\mu}(q)$ for every $w,w' \in W$, we get
\begin{equation*}
\overline{q^{-\langle \rho, \mu \rangle - l(w_0)/2} \sum_{\lambda
\preceq \mu} P_{w_\lambda, w_\mu}(q) \,\, h_\lambda} = q^{-\langle
\rho, \mu \rangle - l(w_0)/2} \sum_{\lambda \preceq \mu} P_{w_\lambda
,w_\mu}(q) \,\, h_\lambda.
\end{equation*} Applying the Satake isomorphism to both sides and using Lemma
\ref{KLcommuteswithSat}, we have
\begin{equation*} q^{\langle \rho, \mu \rangle} \sum_{\lambda \preceq \mu}
P_{w_\lambda, w_\mu}(q^{-1}) \,\, \overline{h_\lambda^\vee} = q^{-\langle
\rho, \mu \rangle} \sum_{\lambda \preceq \mu} P_{w_\lambda ,w_\mu}(q) \,\,
h_\lambda^\vee.
\end{equation*} By (\ref{v-alt.macdonald}),  this gives \vspace{.2in}
\begin{eqnarray*}
\sum_{\lambda \preceq \mu} q^{\langle \rho, \mu - \lambda
\rangle}P_{w_\lambda, w_\mu}(q^{-1}) W_\lambda(q)^{-1} \sum_{w \in W}
t_{w\lambda} \prod_{\alpha > 0} \frac{1 - q \, t_{-w\alpha^\vee}}{1 -
t_{-w\alpha^\vee}}
\hspace{.2in} \,\,\,\,\,\,\,\,\,\,\,\,\,\,\,\,\,\,\,\,\,\, \\ ~ ~ = ~ ~
\sum_{\lambda \preceq \mu} q^{-\langle \rho, \mu -
\lambda \rangle}P_{w_\lambda, w_\mu} (q) W_\lambda (q^{-1})^{-1}
\sum_{w \in W} t_{w\lambda} \prod_{\alpha > 0}
\frac{1-q^{-1}t_{-w\alpha^\vee}}{1 - t_{-w\alpha^\vee}}.
\end{eqnarray*}

Now ${\rm deg}\, P_{w_\lambda, w_\mu}(q) \leq \langle \rho, \mu -
\lambda \rangle - 1/2$ if $\lambda < \mu$ and so Lemma \ref{comb}
 below implies that the right hand side is a polynomial in $q^{-1}$ (with
coefficients in ${\mathbb Z}[X_*]^W$).  Similarly, the left hand side is a
polynomial in $q$.  The result now follows since the constant terms are
equal to
\begin{equation*}
\sum_{w \in W} t_{w\mu} \prod_{\alpha > 0} (1 - t_{-w\alpha^\vee})^{-1} =
E_\mu.
\end{equation*}
\end{proof}

\begin{lemma} \label{comb} We have
\begin{equation*} W_\lambda (q^{-1})^{-1} \sum_{w \in W} t_{w\lambda}
\prod_{\alpha > 0} \frac{1 - q^{-1}t_{-w\alpha^\vee}}{1 - t_{-w\alpha^\vee}}
\in {\mathbb Z}[q^{-1}][X_*]^W.
\end{equation*}
\end{lemma}

\begin{proof} It is obvious that the expression belongs to ${\mathbb
Z}[[q^{-1}]][X_*]^W$, so it is enough by (\ref{v-alt.macdonald}) to show
that $h^\vee_\lambda$ belongs to ${\mathcal R}$.  But this follows from
(\ref{761aaa}).
\end{proof}

Taking $q=1$ we immediately recover Theorem 6.1 of
\cite{lusztig83}:

\begin{theorem} [Lusztig] For any dominant coweight $\mu$,
\begin{equation*} E_\mu = \sum_{\lambda \preceq \mu} P_{w_\lambda, w_\mu}(1)
\,\, (\sum_{w \in W/W_\lambda} t_{w\lambda}).
\end{equation*}
\end{theorem}

\bigskip
\bigskip

\obeylines
University of Maryland
Mathematics Department
College Park, MD 20742-4015
tjh@math.umd.edu

\vspace{.2in}

\obeylines
University of Chicago
Department of Mathematics
5734 S. University Ave.
Chicago, IL 60637 
kottwitz@math.uchicago.edu

\vspace{.2in}

\obeylines

The Institute of Mathematical Sciences
CIT campus Taramani
Chennai 600113, India
amri@imsc.res.in

\end{document}